\documentclass[a4paper,10pt]{article}

\usepackage[a4paper,left=3cm,right=3cm,top=2.5cm,bottom=2.5cm]{geometry}
\usepackage{hyperref}
\usepackage{subfig}
\usepackage{mathtools}
\usepackage{multicol}
\usepackage{comment}
\usepackage{booktabs}
\usepackage{amssymb}
\usepackage{comment}
\usepackage{url}
\usepackage{bm}

\usepackage{pdfsync}
\usepackage{float}
\usepackage{tabularx}
\usepackage{enumerate}
\usepackage{array}
\usepackage{xspace}
\usepackage{siunitx}
\usepackage{authblk}

\usepackage{amsthm}
\usepackage{amsmath}
\usepackage{stmaryrd}
\usepackage[ruled,vlined,linesnumbered]{algorithm2e}
\usepackage{graphicx}
\usepackage[hypcap=true]{subcaption}
\usepackage{booktabs}
\usepackage{cleveref}
\usepackage{dsfont} % for formatting of the indicator function \mathds{1}

\usepackage[normalem]{ulem}

\usepackage{tikz-cd}

\makeatletter
\renewcommand*\env@matrix[1][\arraystretch]{%
  \edef\arraystretch{#1}%
  \hskip -\arraycolsep
  \let\@ifnextchar\new@ifnextchar
  \array{*\c@MaxMatrixCols c}}
\makeatother

\begin{document}

\title{Non-intrusive reduced-order modeling for dynamical systems with spatially localized features}

\author[1]{Leonidas Gkimisis\footnote{gkimisis@mpi-magdeburg.mpg.de}}
\author[2]{Nicole Aretz\footnote{nicole.aretz@austin.utexas.edu}}
\author[3]{Marco Tezzele\footnote{marco.tezzele@emory.edu}}
\author[4]{Thomas Richter\footnote{thomas.richter@ovgu.de}}
\author[1,4]{Peter Benner\footnote{benner@mpi-magdeburg.mpg.de}}
\author[2]{Karen E. Willcox\footnote{kwillcox@oden.utexas.edu}}

\affil[1]{Computational Methods in Systems and Control Theory (CSC), Max Planck Institute for Dynamics of Complex Technical Systems, Sandtorstra{\ss}e 1, 39106 Magdeburg, Germany}

\affil[2]{Oden Institute for Computational Engineering and Sciences, University of Texas
at Austin, Austin, 78712, TX, United States}

\affil[3]{Department of Mathematics, Emory University, Atlanta, 30322, GA, United States}

\affil[4]{Institute for Analysis and Numerics, Otto-von-Guericke-Universit\"at Magdeburg, Universit\"atsplatz 2, 39106 Magdeburg, Germany}

\maketitle

%\maketitle
\begin{abstract}
%\MT{1. Tell the reader what problem you are addressing
%2. Clearly state the contribution of the paper, whether it be a methodological advance or new application capability or new scientific insight.
%3. Give a sense of the evidence that supports the contribution, e.g., the application to which a new method was applied and what was achieved to advance the state of the art.}

This work presents a non-intrusive reduced-order modeling framework for dynamical systems with spatially localized features characterized by slow singular value decay. The proposed approach builds upon two existing methodologies for reduced and full-order non-intrusive modeling, namely Operator Inference (OpInf) and sparse Full-Order Model (sFOM) inference. We decompose the domain into two complementary subdomains that exhibit fast and slow singular value decay. The dynamics of the subdomain exhibiting slow singular value decay are learned with sFOM while the dynamics with intrinsically low dimensionality on the complementary subdomain are learned with OpInf. The resulting, coupled OpInf-sFOM formulation leverages the computational efficiency of OpInf and the high resolution of sFOM, and thus enables fast non-intrusive predictions for conditions beyond those sampled in the training data set. A novel regularization technique with a closed-form solution based on the Gershgorin disk theorem is introduced to promote stable sFOM and OpInf models. We also provide a data-driven indicator for subdomain selection and ensure solution smoothness over the interface via a post-processing interpolation step. We evaluate the efficiency of the approach in terms of offline and online speedup through a quantitative, parametric computational cost analysis. We demonstrate the coupled OpInf-sFOM formulation for two test cases: a one-dimensional Burgers' model for which accurate predictions beyond the span of the training snapshots are presented, and a two-dimensional parametric model for the Pine Island Glacier ice thickness dynamics, for which the OpInf-sFOM model achieves an average prediction error on the order of $1 \%$ with an online speedup factor of approximately $8\times$ compared to the numerical simulation.
\end{abstract}

\section{Introduction}
\label{sec:intro}
Physics-informed non-intrusive modeling leverages partial a priori knowledge about the partial differential equations (PDEs) governing the dynamics of a system to infer models from numerical or experimental data. In this paper, we use domain decomposition to couple two methods for physics-informed non-intrusive modeling that act on both the full and reduced dimension. Thus, we enable accurate predictions on systems with spatially localized features at a reduced computational cost.

Model reduction proposes mathematical frameworks that reduce the computational complexity of numerical simulations while maintaining high accuracy~\cite{morReview, Benner2021}. When a full-order model (FOM), derived via the discretization of governing PDEs, is explicitly available from the numerical solver, so-called ``intrusive'' model reduction methods can be employed~\cite{Benner2015}. However, the numerical simulation of complex or large-scale engineering systems is often performed by legacy code or commercial software. The underlying FOM of the dynamical system at hand can thus be inaccessible to the end-user; simulation data are available, while the governing equations for the system dynamics might be only partially known~\cite{mendez2023datadriven, Benner2021, Karniadakis2021,MG23}.

Physics-informed non-intrusive modeling enables the inference of predictive dynamical models for such cases by leveraging a priori knowledge of the system's governing equations in combination with training data from simulations. These methods offer the capability to make accurate predictions of the underlying system dynamics for different inputs, system parameters, boundary or initial conditions, without accessing the underlying FOM operators. Depending on the availability of simulation training data and the extent to which the underlying governing equations are known, different approaches to physics-informed non-intrusive modeling have been presented over the last two decades such as~\cite{SCHMID_2010, peherstorfer2016data, Raissi2019, karachalios2024datadriven, MG23, SCHUMANN2023, schaeffer2017learning, Brunton2016}. Assuming access to state data, there are projection-based methods that directly infer a reduced-order model (ROM)~\cite{SCHMID_2010, peherstorfer2016data}, and also methods that predict the evolution of the system at the full-order level~\cite{Schaeffer2018, SCHUMANN2023, Raissi2019}. In particular, the methods of Operator Inference (OpInf) \cite{peherstorfer2016data, opinfrev} and sparse FOM (sFOM) inference \cite{Schumann2022, SCHUMANN2023} use state data to infer the system operators in a reduced and the full order accordingly, with the additional assumption of a physics-informed structure of the governing equations.

Projection-based methods offer high online computational efficiency by inferring non-intrusive models directly at a reduced dimension, for example~\cite{pidmd, peherstorfer2016data, XIAO2015522}. Such methods have produced successful predictions for a variety of large-scale applications~\cite{mcquarrie2021data, farcas2023improving, Schmid2011}, while stability of the inferred non-intrusive ROMs has also been addressed~\cite{SAWANT2023115836, Kaptanoglu2021, goyal2023guaranteed}. Projection-based methods are based on the assumption of an intrinsic low dimensionality of the solution manifold, which makes their application to systems with slow singular value decay challenging. Homogeneous systems~\cite{Holmes_Lumley_Berkooz_1996}, transport-dominated systems~\cite{Peherstorfer2020, Rim2023}, or systems driven by the effect of an input or a nonlinearity~\cite{Ohlberger2013} can exhibit slow singular value decay. Several methods have been proposed to tackle this issue focusing on transport-dominated phenomena, via adaptive sampling~\cite{Peherstorfer2020}, coordinate transformations~\cite{Ohlberger2013, Rim2023}, or nonlinear manifold projection~\cite{Geelen2023, goyal2021learning}; however, making predictions for problems with slow singular value decay remains a challenging task because large reduced orders are required for accurately projecting the training data, but also because the system dynamics can evolve beyond the span of the training snapshots.

Analogously to ROMs, studies have investigated the inference of the underlying discretized FOM operators~\cite{Schumann2022, Schaeffer2018, Maddu2023, Gkimisis2024}, with a particular focus on enabling reliable predictions through stability conditions~\cite{SCHUMANN2023, prakash2024datadriven}.
Performing the inference task at the full-order level weakens the dependence on the training data compared to projection-based methods and enables predictions beyond the span of the training snapshot matrix~\cite{SCHUMANN2023, prakash2024datadriven}. This aspect is valuable for the accurate prediction of system dynamics that do not evolve on a low-dimensional manifold  -- a main restriction for ROMs -- but comes at the high computational cost of inference and simulation at the full-order level. Hence, such methods are not easily scalable to large-scale applications.

For many large-scale applications, the dynamics that drive slow singular value decay are spatially localized and thus comprise a family of target applications for domain decomposition methods~\cite{Lucia2001, Buffoni2009, Bergmann2018, GASTALDI1990347}. In intrusive model reduction,~\cite{Lucia2001} presented one of the first efforts for combining ROMs with a FOM for a quasi-1D nozzle flow with a parametrically varying shock position. Several subsequent works have extended the idea of ROM-FOM coupling in different directions. Namely,~\cite{Lucia2003} focused on a spatially coupled ROM-FOM for a 2D hypersonic flow over a blunt body, while~\cite{Sun2008} employed balanced truncation to reduce systems with spatially localized nonlinearities. The authors in~\cite{Buffoni2009} and~\cite{Bergmann2018} proposed iterative and explicit schemes, respectively, for the spatial coupling of a FOM with a Galerkin-free ROM. More recently, approaches for both partitioned, explicit ROM-FOM coupling~\cite{sandia_ROM_FOM_intrusive, DECASTRO2023116398}, as well as implicit coupling via the Schwarz alternating method~\cite{barnett2022schwarz} have been proposed. An optimization-based domain decomposition strategy for incompressible Navier-Stokes was presented in~\cite{Prusak2023} while the authors of~\cite{ddmavris, DIAZ2024} focused on intrusive model reduction with nonlinear manifolds to capture spatially localized dynamics~\cite{huynh2013static, eftang2013port}.

Recent works considered the application of domain decomposition techniques also to non-intrusive MOR. The authors in~\cite{Xiao1, Xiao2} used localized POD bases, interpolated in time via radial basis functions and Gaussian process regression, respectively. Spatially localized OpInf models were considered in~\cite{farcas2023improving} to improve prediction accuracy and reduce the computational requirements for a three-dimensional unsteady rotating detonation rocket engine simulation with $75$ million degrees of freedom. In a similar direction, ~\cite{moore2024dd_opinf} employed domain decomposition for coupling local OpInf models to each other and to first-principle FOMs, utilizing the Schwarz alternating method. The authors in~\cite{ivagnes2024enhancing} proposed a space-dependent aggregation of non-intrusive model predictions comparing different methodologies in aeronautical applications. In the examined contexts, spatially localized ROMs provide advantages over global ROMs although they remain bound to the limitations of projection-based non-intrusive modeling. 

\begin{figure}[!htb]
\centering
  \includegraphics[width=1.\textwidth]{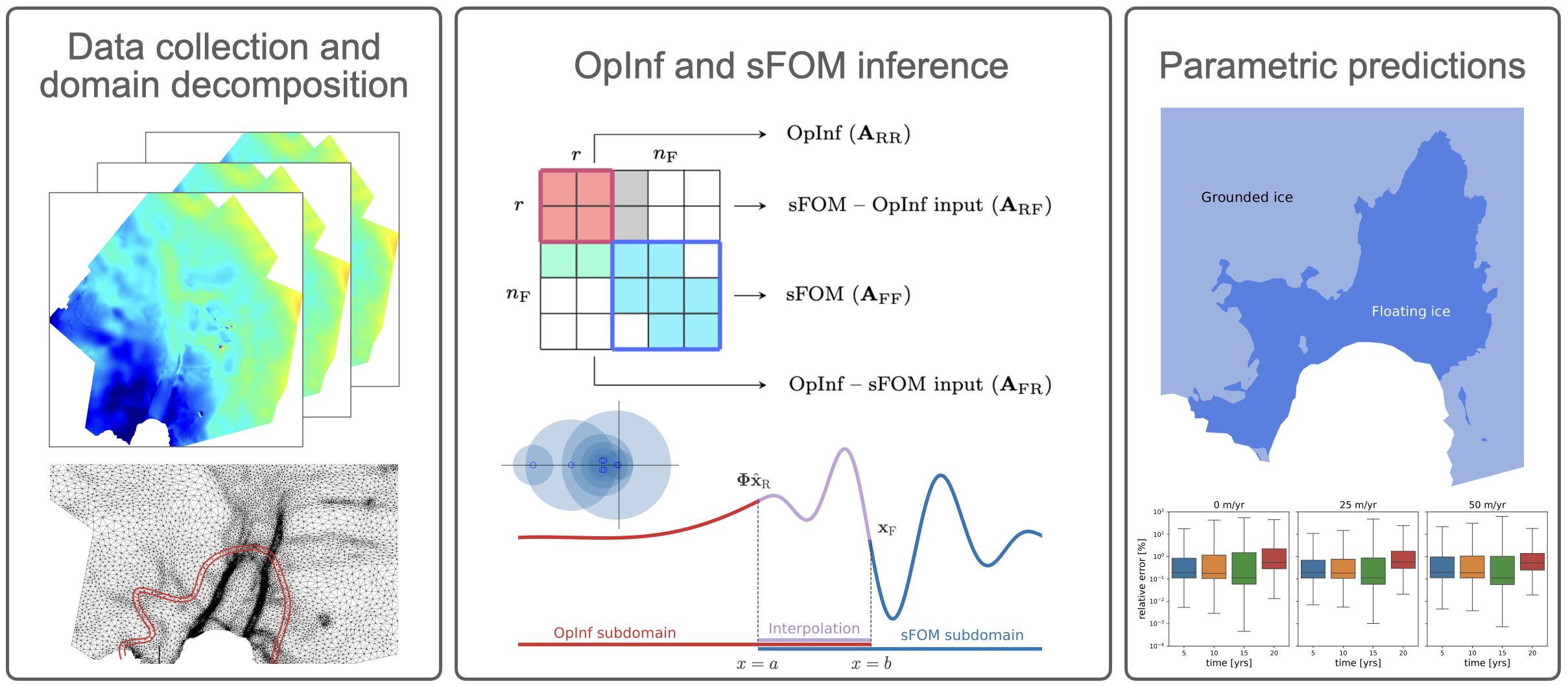}
\caption{Coupled OpInf-sFOM overview. Left panel: we collect state snapshot data and decompose the spatial domain into an OpInf and an sFOM subdomain. Central panel: we perform the coupled OpInf and sFOM inference using a novel stability-promoting regularization and postprocess the interface solution. Right panel: we simulate the coupled OpInf-sFOM and evaluate the predictions.}
\label{fig:paper_scheme}
\end{figure}

In this paper, we build upon the property of spatially localized dynamics driving slow singular value decay in the spirit of~\cite{Lucia2001} to leverage the advantages of a non-intrusive physics-informed ROM and a non-intrusive physics-informed FOM in a domain decomposition framework. An overview of the proposed formulation is given in \Cref{fig:paper_scheme}. For the subdomain where the dynamics evolve in an intrinsically low-dimensional space, we employ the OpInf method~\cite{peherstorfer2016data} to learn a physics-informed non-intrusive ROM, thus achieving high computational efficiency and accuracy. For subdomains where a reduced basis would be too restrictive, we infer a sparse physics-informed non-intrusive FOM via the sFOM approach~\cite{SCHUMANN2023, Gkimisis2024}. The resulting, coupled OpInf-sFOM  extends the prediction capabilities of projection-based non-intrusive ROMs while still achieving a proportional computational speedup with respect to the underlying FOM. To improve the robustness of both models, we introduce a closed-form stability-promoting regularization strategy based on the Gershgorin disk theorem. Moreover, we evaluate the ROM and FOM subdomains selection via the gap in the singular value decay of the corresponding snapshot matrices and achieve solution smoothness across the OpInf-sFOM interface by a post-processing interpolation step. Since the sFOM inference is made at the full-order level, the computational cost is crucial for the scalability of the proposed methodology. We provide quantitative parametric analyses for the offline and online speedup offered by the OpInf-sFOM formulation to a priori estimate the method's computational efficiency. We demonstrate the potential of the OpInf-sFOM approach by making predictions on two numerical test cases. The first is a one-dimensional Burgers' equation example, where the coupled OpInf-sFOM enables the accurate prediction of a propagating wave. The second corresponds to a simplified, 2D parametric model for the Pine Island Glacier ice thickness dynamics. In this case, the coupled parametric OpInf-sFOM achieves accurate ice thickness predictions for a range of input parameters over a period of 20 years, using 15 years of simulation data for the two extreme parameter values for training.

The remainder of this work is organized as follows. \Cref{sec:methods} presents the OpInf and sFOM methods. \Cref{sec:dd} discusses the coupled OpInf-sFOM formulation, solution smoothness across the subdomains interface, and the stability-promoting Gershgorin regularization. An analysis of the offline and online computational cost of the proposed methodology is given in \Cref{sec:cost}. The two numerical test cases along with a discussion on the prediction accuracy are presented in \Cref{sec:results}. Finally, in \Cref{sec:conclusion} we give concluding remarks and potential future research directions.

\section{Theoretical background: Inferring reduced and full-order models}
\label{sec:methods}
The developed methodology comprises a spatial coupling between two non-intrusive modeling methods. On spatial subdomains where the system dynamics can be well approximated within a low-dimensional manifold, we employ OpInf~\cite{peherstorfer2016data}. On spatial subdomains where the projection to a low-dimensional manifold is not sufficiently accurate, such as transport-dominated regions, we use sFOM inference~\cite{Schumann2022, SCHUMANN2023}. We briefly present these two physics-based, non-intrusive methods.

\subsection{State access and physics-informed governing equations}
%{Physics-based non-intrusive modeling}

Both methods build upon two common pillars. First, access to state data is required. We denote the state vector by $\mathbf{x}(t) \in \mathbb{R}^n$ and the input vector by $\mathbf{u}(t) \in \mathbb{R}^k$, where $n$ is the state dimension and $k$ is the input vector dimension. Given data for $\mathbf{x}(t)$ and $\mathbf{u}(t)$ for timesteps $\{ t_1, ..., t_{n_T}\}$, we construct a state snapshot matrix $\mathbf{X} \in \mathbb{R}^{n \times n_T}$ and an input snapshot matrix $\mathbf{U}$ as
\begin{equation}
\label{snapshots}
    \mathbf{X}=\begin{bmatrix} 
    \mathbf{x}(t_1) & \dots & \mathbf{x}(t_{n_T}) 
    \end{bmatrix}, \qquad
    \mathbf{U}=\begin{bmatrix}
    \mathbf{u}(t_1) & \dots & \mathbf{u}(t_{n_T}) 
    \end{bmatrix}. 
\end{equation}
Second, both approaches fall in the category of physics-informed non-intrusive modeling, that is they rely on a priori knowledge of physical laws that govern, exactly or approximately, the system dynamics. This defines a mathematical structure for the underlying discretized model, the dynamics of which we aim to infer. As a prototype, we showcase a quadratic system structure
\begin{equation}
\label{fom_structure}
\frac{ \textrm{d}\mathbf{x}}{\textrm{dt}} =\mathbf{A} \mathbf{x}(t) + \mathbf{H} ( \mathbf{x}(t) \otimes \mathbf{x}(t)) +  \mathbf{B} \mathbf{u}(t) + \mathbf{c}.
\end{equation}
Although the model operators $\mathbf{A} \in \mathbb{R}^{n \times n}, \mathbf{H} \in \mathbb{R}^{n \times n^2}, \mathbf{B} \in \mathbb{R}^{n \times k}, \mathbf{c} \in \mathbb{R}^n$ are unknown, the structure of the discretized dynamical model is fixed. The two methods are presented in the following subsections, illustrating the two different perspectives building upon this common starting point.

\subsection{Operator Inference}
\label{OpInf}

OpInf infers a non-intrusive ROM that encodes the dynamics of a system, given state and input snapshots as well as the governing equations structure~\cite{peherstorfer2016data, mcquarrie2021data}.

Starting from~\eqref{snapshots}, we construct a low-dimensional, linear basis by performing the singular value decomposition (SVD) of the snapshot matrix $\mathbf{X}$, that is
\begin{equation}
\label{svd}
    \mathbf{X} = \mathbf{\Phi} \mathbf{\Sigma} \mathbf{\Psi}^\top,
\end{equation}
where $\mathbf{\Phi} \in \mathbb{R}^{n \times n}$ is the orthonormal matrix formed by the left singular vectors, $\mathbf{\Sigma} \in \mathbb{R}^{n \times n}$ is the diagonal matrix of the singular values in non-increasing order, and $\mathbf{\Psi} \in \mathbb{R}^{n \times n_T}$ is the orthonormal matrix formed by the right singular vectors. Depending on the singular value decay of $\mathbf{X}$, we truncate the first $r \ll n$ left singular vectors in $\mathbf{\Phi}$. The matrix of the first $r$ columns of $\mathbf{\Phi}$, $\mathbf{V} \in \mathbb{R}^{n \times r}$, spans a low-dimensional subspace onto which $\mathbf{x}(t)$ can be projected, leading to a reduced snapshot matrix $\widehat{\mathbf{X}}=  \mathbf{V}^\top \mathbf{X}$. This orthonormal matrix is optimal with respect to the $L_2$ error committed by projection $\mathbf{V} \mathbf{V}^\top \mathbf{x}(t)$~\cite{strang2021introduction}, however, other strategies for constructing low-dimensional bases have been proposed, for example quadratic manifolds~\cite{Geelen2023} or autoencoders~\cite{goyal2021learning}. In the following, we refer to the matrix $\mathbf{V}$ as the reduced basis.

The Galerkin projection of the FOM \eqref{fom_structure} using the reduced basis $\mathbf{V}$ leads to the ROM structure

\begin{equation}
\label{rom_structure}
\frac{\textrm{d}\widehat{\mathbf{x}}}{\textrm{dt}}=\widehat{\mathbf{A}} \widehat{\mathbf{x}}(t) + \widehat{\mathbf{H}} ( \widehat{\mathbf{x}}(t) \otimes \widehat{\mathbf{x}}(t)) +  \widehat{\mathbf{B}} \mathbf{u}(t) + \widehat{\mathbf{c}},
\end{equation}

\noindent with $\widehat{\mathbf{x}}(t) \in \mathbb{R}^r$. The full-order state $\mathbf{x}(t)$ can then be approximated as $\mathbf{x}(t) \approx \mathbf{V} \widehat{\mathbf{x}}(t)$.

Using the projected state and input snapshot data, system operators $\widehat{\mathbf{A}}\in \mathbb{R}^{r \times r}, \widehat{\mathbf{H}}\in \mathbb{R}^{r \times r^2}, \widehat{\mathbf{B}} \in \mathbb{R}^{r \times k}, \widehat{\mathbf{c}} \in \mathbb{R}^r$, are then inferred by solving the least-squares (LS) problem

\begin{equation}
\label{opinf}
\min _{\bm{\mathcal{O}}}{\left\|\bm{\mathcal{O}} \bm{\mathcal{D}}  - \frac{\textrm{d}\widehat{\mathbf{X}}}{\textrm{dt}} \right\|^2_{F}},
\end{equation}
with
\begin{equation}
\label{opinf_operators}
\bm{\mathcal{O}}=\begin{bmatrix}
\widehat{\mathbf{A}}, \; \widehat{\mathbf{H}}, \; \widehat{\mathbf{B}}, \; \widehat{\mathbf{c}} \end{bmatrix}, \qquad 
\bm{\mathcal{D}}=\begin{bmatrix}
\widehat{\mathbf{X}} \\ \widehat{\mathbf{X}} \odot \widehat{\mathbf{X}}\\ \mathbf{U} \\ \widehat{\mathds{1}} \end{bmatrix},
\end{equation}

\noindent where $\odot$ signifies the Khatri-Rao product (columnwise Kronecker product, for each projected snapshot $\widehat{x} (t_i)=V^T  x(t_i)$) and we assume that the state time derivative data $\frac{\textrm{d}\widehat{\mathbf{X}}}{\textrm{dt}}$ are available or can be accurately computed via a chosen numerical scheme from the state data.

We here emphasize that equation~\eqref{opinf} is solved in the reduced dimension $r$ and admits a closed-form solution, making the method computationally efficient, even for large-scale applications~\cite{mcquarrie2021data, farcas2023improving}.  

The inference of $\bm{\mathcal{O}}$ in~\eqref{opinf} allows us to use \eqref{rom_structure} for predicting the system dynamics for initial conditions or input values other than those used to produce the training data in~\eqref{snapshots}. However, all predictions made by inferred models within a projection-based framework are by definition bound to the basis $\mathbf{V}$ and thus the span of the training snapshots. This fact can pose a limitation for accurate ROM predictions of systems with slow singular value decay, such as transport-dominated systems \cite{Peherstorfer2020}.

\subsection{Sparse Full-Order Model Inference}
\label{sFOM}

The sFOM method \cite{SCHUMANN2023} infers the numerical stencil coefficients for each degree of freedom (DOF) in $\mathbf{x}(t)$ for systems with operators that act locally \cite{pidmd,Schumann2022}, given available data in~\eqref{snapshots}, while assuming a physics-based structure from~\eqref{fom_structure} and an available mesh. This leads to the inference of an adjacency-based sparse FOM which interprets the system dynamics \cite{Gkimisis2023, Gkimisis2024}.

To derive the sFOM learning problem, we focus on the governing equation for each entry $\mathbf{x}_i(t)$ of $\mathbf{x}(t)$ individually for $i=1, \dots, n$:
Restricting \eqref{fom_structure} to the $i$-th row, we obtain
\begin{equation}\label{eq:sfom:derivation:1}
    \frac{\textrm{d}\mathbf{x}_i}{\textrm{dt}}(t) 
    % &= \mathbf{A}_{i,:} \mathbf{x}(t) + \mathbf{H}_{i,:} ( \mathbf{x}(t) \otimes \mathbf{x}(t)) +  \mathbf{B}_{i,:} \mathbf{u}(t) + \mathbf{c}_i \\
    = \sum_{j=1}^{n} \mathbf{A}_{i,j} \mathbf{x}_j(t) + \sum_{j=1}^{n^2} \mathbf{H}_{i,j} ( \mathbf{x}(t) \otimes \mathbf{x}(t))_{j} + \sum_{j=1}^{k}  \mathbf{B}_{i,j} \mathbf{u}_j(t) + \mathbf{c}_i.
\end{equation}
Assuming an adjacency-based sparsity structure of the full-order operators, $\mathbf{A}_{i,j}$ can be non-zero only if the $j$-th DOF $\mathbf{x}_j$ is identical or geometrically adjacent to $\mathbf{x}_i$.
We denote this set of indices by $Q_i \subset \{1,\dots,n\}$.
By the same argument, $\mathbf{H}_{i,j}$ is non-zero only for indices $j$ operating on interactions between DOFs in $Q_i$.
We denote this set of indices by $E_i \subset \{1,\dots,n^2\}$.
The sparsity structure of $\mathbf{B}$ depends on the input $\mathbf{u}(t) \in \mathbb{R}^k$;
% Typically, if $k$ is small, e.g., when $\mathbf{u}(t)$ is a scaling factor, $\mathbf{B}$ is dense but skinny;
% if $k \approx n$ is of a similar magnitude as the state dimension $n$, e.g., when $\mathbf{u}(t)$ is a field, then $\mathbf{B}$ is sparse.
to account for both sparse and dense matrices $\mathbf{B}$, we denote the set of non-zero entries of $\mathbf{B}_{i,:}$ by $L_i \subset \{1,\dots,k\}$.
%, stressing that equality is permitted.
With these definitions of the index sets $Q_i, E_i, L_i$ 
% denoting the only entries of $\mathbf{A}_{i,:}, \mathbf{H}_{i,:}, \mathbf{B}_{i,:}$ that can be non-zero,
and with $\mathbf{A}_{i,Q_i} \in \mathbb{R}^{1\times|Q_i|}, \mathbf{H}_{i,E_i} \in \mathbb{R}^{1\times|E_i|}, \mathbf{B}_{i,L_i} \in \mathbb{R}^{1\times|L_i|}$ denoting the restriction of $\mathbf{A}, \mathbf{H}, \mathbf{B}$ onto the $i$-th row and only those columns with index in $Q_i, E_i, L_i$, the summations in \eqref{eq:sfom:derivation:1} collapse to
\begin{equation}\label{eq:sfom:derivation:2}
    \begin{aligned}
    \frac{\textrm{d}\mathbf{x}_i}{\textrm{dt}}(t)
    &= \sum_{j \in Q_i} \mathbf{A}_{i,j} \mathbf{x}_j(t) + \sum_{j \in E_i} \mathbf{H}_{i,j} ( \mathbf{x}(t) \otimes \mathbf{x}(t))_{j} + \sum_{j \in L_i}  \mathbf{B}_{i,j} \mathbf{u}_j(t) + \mathbf{c}_i \\
    % &= \mathbf{A}_{i,Q_i} \mathbf{x}_{Q_i}(t) + \mathbf{H}_{i,E_i} ( \mathbf{x}(t) \otimes \mathbf{x}(t))_{E_i} + \mathbf{B}_{i,L_i} \mathbf{u}_{L_i}(t) + \mathbf{c}_i \\
    &= \mathbf{A}_{i,Q_i} \mathbf{x}_{Q_i}(t) + \mathbf{H}_{i,E_i} ( \mathbf{x}_{Q_i}(t) \otimes \mathbf{x}_{Q_i}(t)) + \mathbf{B}_{i,L_i} \mathbf{u}_{L_i}(t) + \mathbf{c}_i,
\end{aligned}
\end{equation}
where we have used that $( \mathbf{x}(t) \otimes \mathbf{x}(t))_{E_i} = \mathbf{x}_{Q_i}(t) \otimes \mathbf{x}_{Q_i}(t)$ by definition of the set $E_i$.
We further condense \eqref{eq:sfom:derivation:2} to
\begin{equation}
\label{coefs}
    \frac{\textrm{d}\mathbf{x}_i}{\textrm{dt}}(t) =\mathbf{f}^\top_{Q_i}(t) \bm{\beta}_i
\end{equation}
through the auxiliary variables
\begin{equation}
\label{fom_variables_i}
\mathbf{f}_{Q_i}(t) = 
\begin{bmatrix}[1.4]
  \mathbf{x}_{Q_i}(t) \\ \mathbf{x}_{Q_i}(t) \otimes  \mathbf{x}_{Q_i}(t) \\   \mathbf{u}_{L_i}(t) \\ 1
\end{bmatrix},
\qquad
\bm{\beta}_i = 
\begin{bmatrix}[1.4]
  \mathbf{A}_{i, Q_i}^\top \\
  \mathbf{H}_{i, E_i}^\top \\  \mathbf{B}_{i, L_i}^\top \\ \mathbf{c}_i
\end{bmatrix}.
\end{equation}

% \begin{equation}
% \label{fom_variables_i}
% \mathbf{f}_{Q_i}(t) = 
% \begin{bmatrix}[1.4]
%   \mathbf{x}_{Q_i}(t) \\ \mathbf{x}_{Q_i}(t) \odot  \mathbf{x}_{Q_i}(t) \\   \mathbf{u}_{L_i}(t) \\ 1
% \end{bmatrix},
% \qquad
% \bm{\beta}_i = 
% \begin{bmatrix}[1.4]
%   \mathbf{A}_{i, Q_i}^\top \\
%   \mathbf{H}_{i, E_i}^\top \\  \mathbf{B}_{i, L_i}^\top \\ \mathbf{c}_i
% \end{bmatrix},
% \end{equation}

In a non-intrusive setting, the operators $\mathbf{A}, \mathbf{H}, \mathbf{B}, \mathbf{c}$ ~\eqref{fom_structure} are not accessible, and the numerical stencil vector $\bm{\beta}_i$ is thus unknown.
However, formulation~\eqref{fom_variables_i} can be used to infer the operators in~\eqref{fom_structure}, assuming an adjacency-based sparsity pattern, and allows to directly encode Dirichlet boundary condition information. Considering the adjacency-based sparsity pattern for operators $\mathbf{A}, \mathbf{H}, \mathbf{B}, \mathbf{c}$ dramatically reduces the number of entries to be inferred for a FOM. We write the LS problem analogously to~\eqref{opinf} as
\begin{equation}
\label{sfom_ls}
\min_{\bm{\beta}_i}{\left\| \bm{\beta}_i^\top \bm{\mathcal{D}}_i  -  \left. \frac{\textrm{d}\mathbf{X}}{\textrm{dt}} \right\vert_{i} \right\|^2_{2}},
\end{equation}
where
\begin{equation}
\label{data_sfom}
\bm{\mathcal{D}}_i=\begin{bmatrix}
\mathbf{f}_{Q_i}(t_1), \; \dots,\; \mathbf{f}_{Q_i}(t_{n_T})
\end{bmatrix} \; \in \; \mathbb{R}^{k_i \times n_T},
\end{equation}

\noindent is the matrix of auxiliary variable data corresponding to the $k_i$ geometrically adjacent DOFs to DOF $i$, $\left. 
 \frac{\textrm{d}\mathbf{X}}{\textrm{dt}} \right\vert_{i} \in \mathbb{R}^{1 \times n_T}$ is the time derivative vector of DOF $i$ and $\bm{\beta}_i \in \mathbb{R}^{k_i}$ is the vector of unknown numerical coefficients. Equation \eqref{sfom_ls} corresponds to the inference of the $i$-th row of the operators in~\eqref{fom_structure}, given an adjacency-based sparsity pattern. In the special case where the mesh is uniform, data from multiple internal DOFs of the computational domain can be concatenated to~\eqref{data_sfom} and to the  $ \left. \frac{\textrm{d}\mathbf{X}}{\textrm{dt}} \right\vert_{i}$ vector.

The argument of sparsity practically enables inference of the operators at the full-order level, since the computational cost of solving~\eqref{sfom_ls} for every DOF $i$ is independent of the FOM dimension. However, since the total sFOM computational cost scales linearly with the FOM dimension, it is significant for large-scale applications with millions of DOFs. Moreover, the additional requirement of retrieving adjacency information for the system states is not always trivial in cases of proprietary software data. On the other hand, working at the full-order level potentially allows for predictions beyond the space spanned by the training snapshot matrix~\eqref{snapshots}, since the inferred model is not bound to a fixed projection basis. As shown below, this property can be proven useful in systems exhibiting slow singular value decay.

\section{Domain decomposed ROM/FOM formulation}
\label{sec:dd}

The proposed methodology couples a non-intrusive ROM obtained via OpInf with a non-intrusive FOM obtained via sFOM, using domain decomposition. The formulation leverages the predictive capabilities of sFOM for subdomains where the system dynamics exhibit localized features linked to slow singular value decay, with the computational efficiency and accuracy of OpInf for subdomains where dynamics evolve on an intrinsically low-dimensional manifold. In this section, we formulate the coupled problem and discuss aspects of computational cost, ROM and FOM subdomain selection, solution continuity, and stability-promoting regularization.

\subsection{Non-intrusive ROM/FOM formulation via domain decomposition}

We infer a coupled non-intrusive ROM/FOM by splitting the computational domain into two subdomains. \eqref{opinf} is solved in the ROM subdomain and \eqref{sfom_ls} is solved in the FOM subdomain. Additional input terms arise for both OpInf and sFOM from the coupling of the dynamics on the two subdomains. 

\subsubsection{ROM/FOM state vector splitting}

Without loss of generality, we assume that the state vector $\mathbf{x}(t) \in \mathbb{R}^n$ is ordered as

\begin{equation}
\label{split_vector}
  \mathbf{x}(t)=  \begin{bmatrix}
\mathbf{x}_\text{R}(t)\\ \mathbf{x}_\text{F}(t) \end{bmatrix}, 
\end{equation}

\noindent where $\mathbf{x}_\text{R}(t) \in \mathbb{R}^{n_\text{R}}$ is the vector of DOFs in the subdomain attributed to the ROM and $\mathbf{x}_\text{F}(t)  \in \mathbb{R}^{n_\text{F}}$ is the vector of DOFs in the FOM subdomain. Snapshot matrices $\mathbf{X}_\text{R} \in \mathbb{R}^{n_R \times n_T}$ and $\mathbf{X}_\text{F} \in \mathbb{R}^{n_F \times n_T}$ are assembled in accordance to \eqref{snapshots}.

We then compute a reduced basis for the ROM subdomain from the snapshot matrix of $\mathbf{x}_\text{R}(t)$. Following the exposition in \Cref{OpInf}, we denote the reduced basis for the ROM subdomain by $\mathbf{V}_\text{R}  \in \mathbb{R}^{n_\text{R} \times r}$ and the corresponding projected snapshot matrix by $\widehat{\mathbf{X}}_\text{R}$, such that we can write

\begin{equation}
\label{red_vector}
\mathbf{x}(t) \approx  \begin{bmatrix}
\mathbf{V}_\text{R} \widehat{\mathbf{x}}_\text{R}(t)\\ \mathbf{x}_\text{F}(t) \end{bmatrix}.
\end{equation}

\noindent We then substitute~\eqref{red_vector} into the physics-informed structure for the dynamical evolution of $\mathbf{x}(t)$. This leads to a physics-informed, coupled ROM-FOM structure, which we aim to learn via OpInf and sFOM.

A critical aspect of the decomposition problem in \eqref{red_vector} is the selection of the ROM and FOM subdomains. In this study, we limit ourselves to an a priori selection of the OpInf and sFOM subdomains.  For this selection, we monitor the singular value decay of the data in both subdomains, such that the decay in the sFOM subdomain is slower than that in the OpInf subdomain. This is done by inspecting the training data and identifying regions where phenomena linked to slow singular value decay arise, such as transport-dominated dynamics or localized nonlinearities. These regions are then assigned to the sFOM subdomain. An adaptive update of the ROM and FOM subdomains based on an a posteriori error estimator is left as a future research endeavour.

\subsubsection{Linear autonomous systems}

For ease of exposition, we first consider the model structure in~\eqref{fom_structure} with only a linear term

\begin{equation}
\label{fom_linear}
\frac{\textrm{d}\mathbf{x}}{\textrm{dt}}=\mathbf{A} \mathbf{x}(t),
\end{equation}
    
By employing~\eqref{red_vector} in~\eqref{fom_linear} and considering the Galerkin projection of $\mathbf{x}_\text{R}(t)$, we obtain

\begin{equation}
\begin{dcases}
\begin{aligned}
\label{coupled_structure_linear}
\frac{\textrm{d}\widehat{\mathbf{x}}_\text{R}}{\textrm{dt}} &= \mathbf{A}_{\text{RR}} \widehat{\mathbf{x}}_\text{R}(t) + \mathbf{A}_{\text{RF}} {\mathbf{x}}_\text{F}(t)
\\ %\;\\
\frac{\textrm{d}{\mathbf{x}}_\text{F}}{\textrm{dt}} &= \mathbf{A}_{\text{FF}} {\mathbf{x}}_\text{F}(t) + \mathbf{A}_{\text{FR}} \widehat{\mathbf{x}}_\text{R}(t),
\end{aligned}
\end{dcases}
\end{equation}

\noindent where $\mathbf{A}_{\text{RR}} \in \mathbb{R}^{r \times r}$, $\mathbf{A}_{\text{RF}} \in \mathbb{R}^{r \times n_\text{F}}$, $\mathbf{A}_{\text{FF}} \in \mathbb{R}^{n_\text{F} \times n_\text{F}}$ and $\mathbf{A}_{\text{FR}} \in \mathbb{R}^{n_\text{F} \times r}$. The first letter of the linear matrix subscript denotes the (ROM or FOM) model to which it refers and the second letter denotes the (ROM or FOM) vector with which it is multiplied. We observe that \eqref{coupled_structure_linear} is a two-way coupled system with a linear coupling input term in each subsystem.

\begin{figure}[!htb]
\centering
\begin{tikzpicture}[every node/.style={minimum size=.5cm-\pgflinewidth, outer sep=0pt}]

    %%%%%%%% first matrix
    \node at (-.8, 1.25) {$\mathbf{A} = $};
    % top side
    \node at (1., 3.2) {\footnotesize $n_\text{R}$};
    \node at (2.75, 3.2) {\footnotesize $n_\text{F}$};
    % left side
    \node at (-.3, 2.) {\footnotesize $n_\text{R}$};
    \node at (-.3, .25) {\footnotesize $n_\text{F}$};
    
    % main grid
    \draw[step=0.5cm,color=black, shift={(0,-.5)}] (0,0) grid (3.5,3.5);

    % 7th row
    \node[draw=darkgray, line width=0.12mm, fill=blue!10!cyan!40] at (2.75, -0.25) {};
    \node[draw=darkgray, line width=0.12mm, fill=blue!10!cyan!40] at (3.25, -0.25) {};
    % 6th row
    \node[draw=darkgray, line width=0.12mm, fill=blue!10!cyan!40] at (2.25, 0.25) {};
    \node[draw=darkgray, line width=0.12mm, fill=blue!10!cyan!40] at (2.75, 0.25) {};
    \node[draw=darkgray, line width=0.12mm, fill=blue!10!cyan!40] at (3.25, 0.25) {};
    % 5th row
    \node[draw=darkgray, line width=0.12mm, fill=blue!10!cyan!40] at (1.75, 0.75) {};
    \node[draw=darkgray, line width=0.12mm, fill=blue!10!cyan!40] at (2.25, 0.75) {};
    \node[draw=darkgray, line width=0.12mm, fill=blue!10!cyan!40] at (2.75, 0.75) {};
    % 4th row
    \node[draw=darkgray, line width=0.12mm, fill=blue!10!cyan!40] at (1.25, 1.25) {};
    \node[draw=darkgray, line width=0.12mm, fill=blue!10!cyan!40] at (1.75, 1.25) {};
    \node[draw=darkgray, line width=0.12mm, fill=blue!10!cyan!40] at (2.25, 1.25) {};
    % 3rd row
    \node[draw=darkgray, line width=0.12mm, fill=blue!10!cyan!40] at (0.75, 1.75) {};
    \node[draw=darkgray, line width=0.12mm, fill=blue!10!cyan!40] at (1.25, 1.75) {};
    \node[draw=darkgray, line width=0.12mm, fill=blue!10!cyan!40] at (1.75, 1.75) {};
    % 2nd row
    \node[draw=darkgray, line width=0.12mm, fill=blue!10!cyan!40] at (0.25, 2.25) {};
    \node[draw=darkgray, line width=0.12mm, fill=blue!10!cyan!40] at (0.75, 2.25) {};
    \node[draw=darkgray, line width=0.12mm, fill=blue!10!cyan!40] at (1.25, 2.25) {};
    % 1st row
    \node[draw=darkgray, line width=0.12mm, fill=blue!10!cyan!40] at (0.25, 2.75) {};
    \node[draw=darkgray, line width=0.12mm, fill=blue!10!cyan!40] at (0.75, 2.75) {};

    % grids
    \draw[step=2cm, color=purple!70, line width=.6mm, shift={(0,1)}] (0,0) grid (2,2);
    \draw[step=1.5cm, color=blue!80!cyan!70, line width=.6mm, shift={(2,-.5)}] (0,0) grid (1.5,1.5);

    % arrow old version
    %\node at (4.2, 1.5) {\footnotesize DD};
    %\node at (4.2, 1.0) {\footnotesize ROM};
    %\draw [-stealth] (3.8, 1.25) -- (4.7, 1.25);
    
    % arrow
    \node at (5., 1.9) {\footnotesize Domain};
    \node at (5., 1.5) {\footnotesize decomposition};
    \node at (5., 0.95) {\footnotesize Projection};
    \draw [-stealth, line width=.2mm] (3.8, 1.25) -- (6.2, 1.25);

    %%%%%%%% second matrix
    % \node at (5.3, 1.3) {$\widehat{\mathbf{A}} = $}; % old version
    % top side
    \node at (7.2, 2.7) {\footnotesize $r$};
    \node at (8.45, 2.7) {\footnotesize $n_\text{F}$};
    % left side
    \node at (6.4, 2) {\footnotesize $r$};
    \node at (6.4, .75) {\footnotesize $n_\text{F}$};

    % main grid
    \draw[step=0.5cm, color=black, shift={(6.7,0)}] (0,0) grid (2.5,2.5);

    % 5th row
    \node[draw=darkgray, line width=0.12mm, fill=blue!10!cyan!50] at (8.45, 0.25) {};
    \node[draw=darkgray, line width=0.12mm, fill=blue!10!cyan!50] at (8.95, 0.25) {};
    % 4th row
    \node[draw=darkgray, line width=0.12mm, fill=blue!10!cyan!50] at (7.95, 0.75) {};
    \node[draw=darkgray, line width=0.12mm, fill=blue!10!cyan!50] at (8.45, 0.75) {};
    \node[draw=darkgray, line width=0.12mm, fill=blue!10!cyan!50] at (8.95, 0.75) {};
    % 3rd row
    \node[draw=darkgray, line width=0.12mm, fill=green!40!cyan!40] at (6.95, 1.25) {};
    \node[draw=darkgray, line width=0.12mm, fill=green!40!cyan!40] at (7.45, 1.25) {};
    \node[draw=darkgray, line width=0.12mm, fill=blue!10!cyan!50] at (7.95, 1.25) {};
    \node[draw=darkgray, line width=0.12mm, fill=blue!10!cyan!50] at (8.45, 1.25) {};
    % 2nd row
    \node[draw=darkgray, line width=0.12mm, fill=red!40] at (6.95, 1.75) {};
    \node[draw=darkgray, line width=0.12mm, fill=red!40] at (7.45, 1.75) {};
    \node[draw=darkgray, line width=0.12mm, fill=black!20] at (7.95, 1.75) {};
    % 1st row
    \node[draw=darkgray, line width=0.12mm, fill=red!40] at (6.95, 2.25) {};
    \node[draw=darkgray, line width=0.12mm, fill=red!40] at (7.45, 2.25) {};
    \node[draw=darkgray, line width=0.12mm, fill=black!20] at (7.95, 2.25) {};

    % grids
    \draw[step=1cm, color=purple!70, line width=.6mm, shift={(6.7,1.5)}] (0,0) grid (1,1); % opinf
    \draw[step=1.5cm, color=blue!80!cyan!70, line width=.6mm, shift={(7.7,0)}] (0,0) grid (1.5,1.5); % sfom

    % arrows
    % (10.5, .75)
    \node at (10.5, .75) {\footnotesize $\mathbf{A}_{\text{FF}}$ };
    \draw [-stealth, color=blue!80!cyan!70, line width=.3mm] (9.4, .75) -- (9.85, .75);
    % (sFOM -- OpInf input) (11.5, 2.)
    \node at (10.5, 2.) {\footnotesize $\mathbf{A}_{\text{RF}}$ };
    \draw [-stealth, color=black!50, line width=.3mm] (9.4, 2.) -- (9.85, 2.);
%(10.5, 3.)
    \node at (10.5, 3.) {\footnotesize $\mathbf{A}_{\text{RR}}$ };
    \draw [-stealth, color=purple!70, line width=.3mm] (7.45, 3.) -- (9.85, 3.);
    \draw [color=purple!70, line width=.3mm] (7.45, 2.6) -- (7.45, 3.);
    %(OpInf -- sFOM input) (11.5, -.5)
    \node at (10.5, -.5) {\footnotesize $\mathbf{A}_{\text{FR}}$ };
    \draw [-stealth, color=green!50!cyan!90!blue!70, line width=.3mm] (7.45, -.5) -- (9.85, -.5);
    \draw [color=green!50!cyan!90!blue!70, line width=.3mm] (7.45, -.1) -- (7.45, -.5);

  \end{tikzpicture} 

\caption{Schematic representation of ROM/FOM domain decomposition for linear operator $\mathbf{A}$. $\mathbf{A}$ is assumed to be sparse. Projecting the first $n_\text{R}$ DOFs to a ROM and maintaining $n_\text{F}$ DOFs at the FOM level results in a coupled structure, with a dense $\mathbf{A}_{\text{RR}} \in \mathbb{R}^{r \times r}$, a sparse $\mathbf{A}_{\text{FF}} \in \mathbb{R}^{n_\text{F} \times n_\text{F}}$, and two coupling operators $\mathbf{A}_{\text{FR}} \in \mathbb{R}^{n_\text{F} \times r}$ and $\mathbf{A}_{\text{RF}} \in \mathbb{R}^{r \times n_\text{F}}$. We learn these operators via~\eqref{coupled_inference}.}
    \label{fig:matrix_structures}
\end{figure}
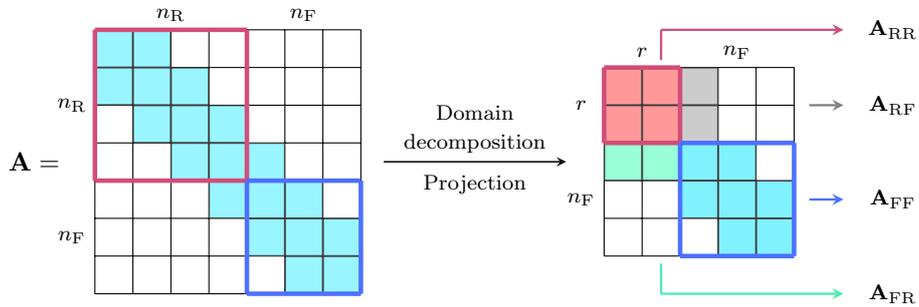

The splitting of the linear matrix $\mathbf{A}$ is sketched in \Cref{fig:matrix_structures}. \Cref{fig:matrix_structures} indicates that the coupling input from the FOM to the ROM originates from the subset $\text{I}$ of FOM DOFs located in the vicinity of the subdomains' interface, $\mathbf{x}_\text{I}(t) \in \mathbb{R}^{n_\text{I}}$, where $n_\text{I} \ll n_\text{F}$. Hence, the term $\mathbf{A}_{\text{RF}}{\mathbf{x}}_\text{F}(t)$ is simplified to $\mathbf{A}_{\text{RI}}{\mathbf{x}}_\text{I}(t)$, where $\mathbf{A}_{\text{RI}} \in \mathbb{R}^{r \times n_\text{I}}$. We employ OpInf for the inference of the $\widehat{\mathbf{x}}_\text{R}(t)$ dynamics and sFOM for the inference of the ${\mathbf{x}}_\text{F}$ dynamics in~\eqref{coupled_structure_linear}. The coupled inference problem then reads as
\begin{equation}
\begin{dcases}
\begin{aligned}
\label{coupled_inference}
&\min _{\bm{\mathcal{O}}_\text{R}}{\left\|\bm{\mathcal{O}}_\text{R} \bm{\mathcal{D}}_\text{R}  - \frac{\textrm{d}\widehat{\mathbf{X}}_\text{R}}{\textrm{dt}} \right\|^2_{F}} \\ 
&\min _{\bm{\beta}_i}{\left\|\bm{\beta}_i^\top {\bm{\mathcal{D}}_{\text{F}_i}} -  \left. \frac{\textrm{d}\mathbf{X}_\text{F}}{\textrm{dt}} \right\vert_{i} \right\|^2_{2}}, \quad \forall i \in \left[ 1, \dots, n_F\right],
\end{aligned}
\end{dcases}
\end{equation}

\noindent with

\begin{equation*}
\label{data_coupled_linear}
\bm{\mathcal{O}}_\text{R}=\begin{bmatrix}[1.4]
{\mathbf{A}_{\text{RR}}}, \;{\mathbf{A}_{\text{RI}}} \end{bmatrix}, 
\quad 
\bm{\mathcal{D}}_\text{R}=\begin{bmatrix}[1.4]
\widehat{\mathbf{X}}_\text{R} \\ {\mathbf{X}}_\text{I} \end{bmatrix},
\quad
\bm{\beta}_i = 
\begin{bmatrix}[1.4]
  \mathbf{A}_{\text{FF}_{i, Q_i}}^\top \\  \mathbf{A}_{\text{FR}_{i, :}}^\top
  %\mathbf{A}_{\text{FF}_{i, Q_i}},\; \delta_i^\text{I} \mathbf{A}_{\text{FR}_{i, :}},\;  \mathbf{B}_{\text{F}_{i, L_i}},\; \mathbf{c}_{\text{F}_i}^\top
  \end{bmatrix},
\quad
\bm{\mathcal{D}}_{\text{F}_i}= 
\begin{bmatrix}[1.4]
  {\mathbf{X}_{\text{F}_{Q_i}}} \\ \delta_i^\text{I} {\widehat{\mathbf{X}}}_\text{R}  
\end{bmatrix} ,
\end{equation*}

\noindent where $\bm{\mathcal{O}}_\text{R} \in \mathbb{R}^{ r \times \left(r +n_\text{I} \right) }$, $\bm{\mathcal{D}}_\text{R} \in \mathbb{R}^{ \left(r +n_\text{I} \right) \times n_t}$, $\bm{\mathcal{\beta}}_i \in \mathbb{R}^{ n_{Q_i} + r }$ and $\bm{\mathcal{D}}_{\text{F}_i} \in \mathbb{R}^{ \left( n_{Q_i} + r \right) \times n_t}$ and the Kronecker delta $\delta_i^\text{I}$ is introduced to signify that the coupling FOM/ROM terms are active only for DOFs $i$ in the vicinity of the interface $i \in \text{I}$.

We hereby mention two potential variations to the formulation in \eqref{coupled_inference}. First, we note that for the OpInf problem (the first line in \eqref{coupled_inference}), an input term with only a subset of the $n_\text{I}$ sFOM DOFs could be considered, to reduce the dimension of $\bm{\mathcal{O}}_\text{R}$. Secondly, the sFOM inference could be performed without explicitly using the OpInf vector for $\mathbf{A}_{\text{FR}} \widehat{\mathbf{x}}_\text{R}(t)$. Instead, $\mathbf{A}_{\text{FR}}$ can be computed in a two-step process: first infer an input term at the full-order level from the data of adjacent DOFs of  $i \in \text{F}$ that belong to the OpInf subdomain and then project it to the ROM basis. The input from the ROM to the FOM is then given by the corresponding rows of $\mathbf{V} \widehat{\mathbf{X}}_\text{R}$. This reduces the size of 
$\bm{\beta}_i$ and bypasses potential scaling issues between $\widehat{\mathbf{X}}_\text{R}$ and $\mathbf{X}_\text{F}$ for the numerical solution of~\eqref{coupled_inference}.

\subsubsection{Nonlinear systems}

For systems with higher order polynomial nonlinearities, the exposition becomes more involved since more coupling terms arise from the ROM-FOM decomposition. Following an analogous approach to~\eqref{coupled_structure_linear}, we split the state vector and the input vector in~\eqref{fom_structure} into ROM and FOM components. The resulting coupled system for a model with quadratic nonlinearity, reads as

\begin{equation}
\begin{dcases}
\begin{aligned}
\label{coupled_structure_quad}
\frac{\textrm{d}\widehat{\mathbf{x}}_\text{R}}{\textrm{dt}} &= \mathbf{A}_{\text{RR}} \widehat{\mathbf{x}}_\text{R}(t) + \mathbf{A}_{\text{RF}} {\mathbf{x}}_\text{F}(t) +
\mathbf{H}_{\text{RRR}} (\widehat{\mathbf{x}}_\text{R}(t)\otimes \widehat{\mathbf{x}}_\text{R}(t)) +
\mathbf{H}_{\text{RFF}} ({\mathbf{x}}_\text{F}(t)\otimes {\mathbf{x}}_\text{F}(t)) \,+ \\ & \quad  +
 \mathbf{H}_{\text{RRF}} (\widehat{\mathbf{x}}_\text{R}(t)\otimes {\mathbf{x}}_\text{F}(t)) +  \mathbf{B}_\text{R} \mathbf{u}_\text{R}(t) + \mathbf{c}_\text{R} ,
\\
\frac{\textrm{d}{\mathbf{x}}_\text{F}}{\textrm{dt}} &= \mathbf{A}_{\text{FF}} {\mathbf{x}}_\text{F}(t) + \mathbf{A}_{\text{FR}} \widehat{\mathbf{x}}_\text{R}(t) +
\mathbf{H}_{\text{FFF}} ({\mathbf{x}}_\text{F}(t)\otimes{\mathbf{x}}_\text{F}(t)) +
\mathbf{H}_{\text{FRR}} (\widehat{\mathbf{x}}_\text{R}(t)\otimes \widehat{\mathbf{x}}_\text{R}(t)) \,+ \\ & \quad + \mathbf{H}_{\text{FFR}} ({\mathbf{x}}_\text{F}(t)\otimes \widehat{\mathbf{x}}_\text{R}(t)) +  \mathbf{B}_\text{F} \mathbf{u}_\text{F}(t) + \mathbf{c}_\text{F}.
\end{aligned}
\end{dcases}
\end{equation}

\noindent As in~\eqref{coupled_structure_linear}, the first letter in the matrix subscript denotes whether the matrix corresponds to the FOM (F) or the ROM (R). For $\mathbf{H}$, the last two letters denote the first and second vector of the Kronecker product, hence $\mathbf{H}_{\text{RRF}}$ and $\mathbf{H}_{\text{FFR}}$ introduce bilinear coupling terms for the ROM and the FOM, accordingly. Equation \eqref{coupled_structure_quad} shows that for a system with quadratic nonlinearity, the coupling terms that arise are linear, bilinear and quadratic. Corresponding coupling terms appear for systems with polynomial nonlinearities of higher order. For the quadratic case, the operators in~\eqref{coupled_inference} are as follows:

\begin{equation*}
\label{data_coupled_quad}
\bm{\mathcal{O}}_\text{R}=\begin{bmatrix}[1.4]
\mathbf{A}_{\text{RR}}^\top \\ \mathbf{A}_{\text{RI}}^\top \\ \mathbf{H}_{\text{RRR}}^\top\\ \mathbf{H}_{\text{RII}}^\top\\ \mathbf{H}_{\text{RRI}}^\top\\ \mathbf{B}_\text{R}^\top \\ \mathbf{c}_\text{R}^\top\end{bmatrix}^\top, 
\quad 
\bm{\mathcal{D}}_\text{R}=\begin{bmatrix}[1.4]
\widehat{\mathbf{X}}_\text{R} \\ {\mathbf{X}}_\text{I} \\ \widehat{\mathbf{X}}_\text{R} \odot \widehat{\mathbf{X}}_\text{R}\\ {\mathbf{X}}_\text{I} \odot {\mathbf{X}}_\text{I}\\\widehat{\mathbf{X}}_\text{R} \odot {\mathbf{X}}_\text{I}\\ \mathbf{U}_\text{R} \\ \widehat{\mathds{1}}_\text{R} \end{bmatrix},
\quad
\bm{\beta}_i = 
\begin{bmatrix}[1.4]
  \mathbf{A}_{\text{FF}_{i, Q_i}}^\top \\ \mathbf{A}_{\text{FR}_{i, :}}^\top \\ \mathbf{H}_{\text{FFF}_{i, E_i}}^\top \\ \mathbf{H}_{\text{FRR}_{i, :}}^\top \\ \mathbf{H}_{\text{FFR}_{i, G_i}}^\top \\ \mathbf{B}_{\text{F}_{i, L_i}}^\top \\ \mathbf{c}_{\text{F}_i}
  \end{bmatrix},
\quad
\bm{\mathcal{D}}_{\text{F}_i}= 
\begin{bmatrix}[1.4]
  {\mathbf{X}_{\text{F}_{Q_i}}} \\ \delta_i^\text{I} {\widehat{\mathbf{X}}}_\text{R} \\ {\mathbf{X}_{\text{F}_{Q_i}}} \odot {\mathbf{X}_{\text{F}_{Q_i}}} \\
  \delta_i^\text{I} {\widehat{\mathbf{X}}}_\text{R} \odot {\widehat{\mathbf{X}}}_\text{R}\\
  \delta_i^\text{I} {\mathbf{X}_{\text{F}_{Q_i}}} \odot {\widehat{\mathbf{X}}}_\text{R} \\ \mathbf{U}_{\text{F}_{L_i}} \\ \mathds{1}_\text{F}
\end{bmatrix}.
\end{equation*}

\noindent Analogously to the linear case \eqref{coupled_structure_linear}, we simplify the OpInf coupling operators, by using the FOM DOFs in the vicinity of the interface, denoted by I. The Kronecker delta $\delta_i^\text{I}$ is introduced for the coupling terms in sFOM inference. Entries $E_i$ correspond to the combinations of DOFs $Q_i$, while $L_i$ includes any geometrically adjacent input DOFS to $i$. Entries $G_i$ of $\mathbf{H}_{\text{FFR}}$ result from the combinations of DOFs $Q_i$ with the $r$ OpInf DOFs. The dimensions of the matrices and vectors in \eqref{data_coupled_quad} are

\begin{align*}
\bm{\mathcal{O}}_\text{R} &\in \mathbb{R}^{ r \times \left(r +n_\text{I} + r^2 + n_\text{I}^2 + r \times n_\text{I} + k_\text{R} + 1 \right) }, \\
\bm{\mathcal{D}}_\text{R} &\in \mathbb{R}^{ \left(r +n_\text{I} + r^2 + n_\text{I}^2 + r \times n_\text{I} + + k_\text{R} + 1 \right) \times n_T}, \\ \bm{\mathcal{\beta}}_i &\in \mathbb{R}^{ \left( n_{Q_i} +  r + n_{E_i} + r^2 + n_{G_i} + n_{L_i} +1 \right)}, \\
\bm{\mathcal{D}}_{\text{F}_i} &\in \mathbb{R}^{ \left( n_{Q_i} + r + n_{E_i} + r^2 + n_{G_i} + n_{L_i} +1 \right) \times n_T}. 
\end{align*}
The non-unique DOFs in the polynomial operator entries can be eliminated due to the Kronecker product permutations. For instance, the unique DOFs of $\mathbf{H}_{\text{RRR}}$ are $r(r+1)/2$ out of the total $r^2$ columns of the matrix.

Finally, we note the possibility of eliminating terms in either the OpInf or the sFOM structure of~\eqref{coupled_structure_quad} based on the system dynamics in either subdomain. For example, in incompressible fluid dynamics problems, the linear viscosity term is mainly active in the flow boundary layer, while it is theoretically negligible in the inviscid far field \cite{Noack2003}.

\subsubsection{ROM/FOM interface solution smoothness}

As presented in \Cref{OpInf}, the OpInf model is bound to the space spanned by basis $\mathbf{V}_\text{R}$, while the inferred sFOM operates at the full-order level. Hence, even after successful inference on both sides, the spatially coupled solution could be insufficiently smooth across the ROM/FOM interface.

\begin{figure}[htb!]
\centering
\includegraphics[width=0.85\textwidth]{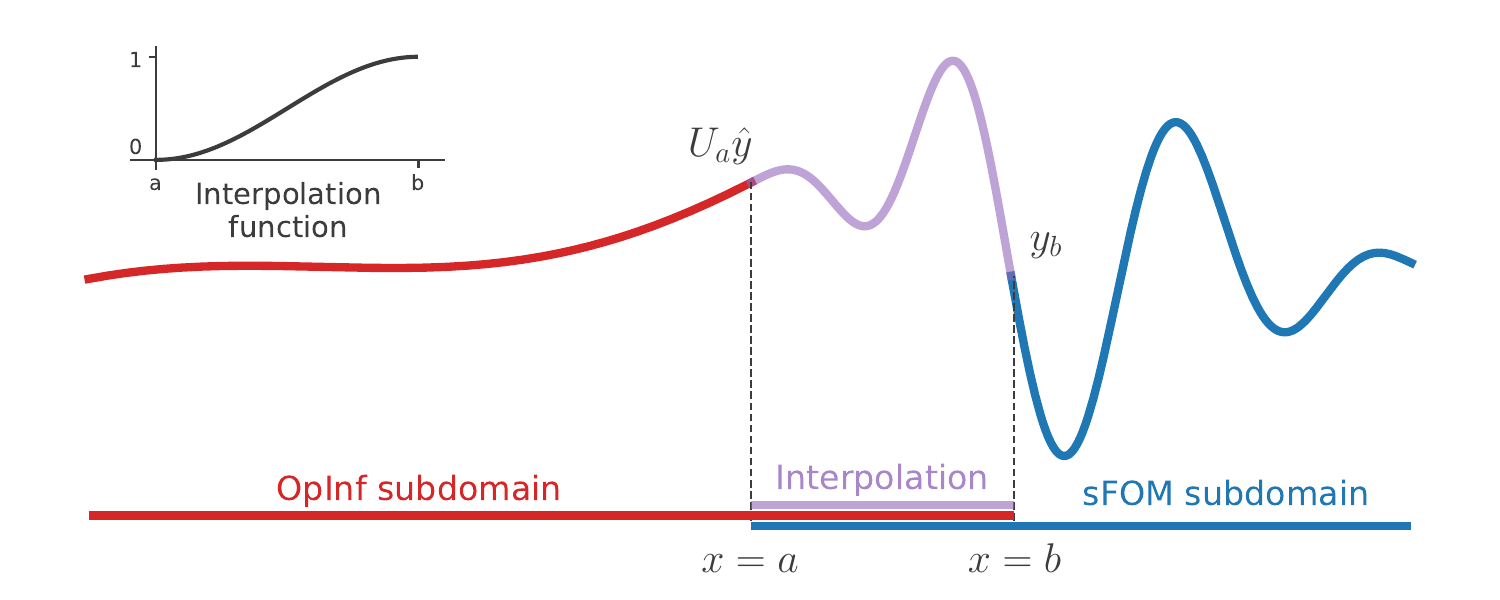}
\caption{Schematic of solution interpolation across the ROM/FOM interface (shown for 1D problems). A transition from the reprojected ROM solution to the FOM solution is ensured by a post-processing interpolation over the overlap region $z \in [a, b]$.}
\label{fig:overlap}
\end{figure}

To alleviate this issue, we select overlapping FOM and ROM subdomains (see \Cref{fig:overlap}) and set the coupling DOFs for either side of~\eqref{coupled_inference} accordingly. As a post-processing step of the coupled OpInf-sFOM predictions, we interpolate the solution in the overlap region (see \Cref{fig:overlap}). The interpolation is performed in a similar fashion to the post-processing interpolation proposed in \cite{farcas2023improving}, where the state at the overlap region between two subdomains is computed as a weighted linear combination of the reconstructed predictions made by the respective non-intrusive models. Alternatively, a bilinear interpolation between the ROM and sFOM values at the interface boundaries, $z=a$, $z=b$ ($z$ denoting the spatial coordinate) can be readily applied to problems with a narrow overlap region without significant loss of accuracy. This has been found sufficient for a smooth overlap solution in the later presented 2D test case.

\subsection{Stability-promoting regularization}

The resulting OpInf~\eqref{sfom_ls} and sFOM~\eqref{opinf} LS problems are often badly conditioned and thus, we complement them with a regularization term~\cite{mcquarrie2021data, farcas2023improving, Gkimisis2023, prakash2024datadriven}. In this section we propose a stability-promoting variant of the $L_2$ regularization strategy based on the Gershgorin disk theorem that has a closed-form solution, for both OpInf and sFOM inference. 

\subsubsection{Gershgorin disks connection to $L_2$ regularization}

The most widely used type of regularization, for both the OpInf~\eqref{opinf} and the sFOM~\eqref{sfom_ls} LS problems, is a Tikhonov or $L_2$ regularization, which penalizes the $L_2$ norm of the operators. For the coupled OpInf-sFOM in~\eqref{coupled_inference} this is written as

\begin{equation}
	\begin{dcases}
		\begin{aligned}
			\label{coupled_reg}
			&\min _{\bm{\mathcal{O}}_\text{R}}{\left\|\bm{\mathcal{O}}_\text{R} \bm{\mathcal{D}}_\text{R}  - \frac{ \textrm{d}\widehat{\mathbf{X}}_\text{R}}{\textrm{dt}} \right\|^2_F + \eta_\text{R}\left\| \bm{\mathcal{O}}_\text{R} \right\|^2_{F} } , \\ 
			&\min _{\bm{\beta}_i}{\left\|\bm{\beta}_i^\top {\bm{\mathcal{D}}_{\text{F}_i}} -  \left. \frac{ \textrm{d}\mathbf{X}_\text{F}}{\textrm{dt}} \right\vert_{i} \right\|^2_{2} + \eta_i \left\|\bm{\beta}_i \right\|^2_2 }, \quad \forall i \in \left[ 1, \dots, n_F\right],
		\end{aligned}
	\end{dcases}
\end{equation}

\noindent where $\eta_\text{R}$ and $\eta_i$ are weighting the regularization term in either LS problem. ~\eqref{coupled_reg} admits a closed-form solution, however the resulting OpInf or sFOM is sensitive to regularization, rendering the selection of $\eta$ values crucial for successful inference~\cite{Gkimisis2023, mcquarrie2021data, benner_flow}. We briefly summarize the selection approach based on an L-curve criterion~\cite{Hansen2000}, where either LS problem in~\eqref{coupled_reg} is solved for different $\eta$ values. The value for which the (either OpInf or sFOM) $L_2$ error in~\eqref{coupled_inference} and the corresponding solution norm are leveraged in terms of a Pareto front selection criterion, qualifies as the optimal $\eta$. The regularization of the LS solution norm is not always uniform. Depending on the system structure in~\eqref{coupled_structure_quad}, the solution entries correspond to different terms of the dynamical system and thus the regularization of their Frobenius norm should be accordingly scaled~\cite{SAWANT2023115836, mcquarrie2021data}.

From a numerical linear algebra viewpoint, the $L_2$ regularization term reduces the condition number of the data matrix in~\eqref{coupled_reg} \cite{strang2021introduction}. Conversely, from a dynamical systems viewpoint, the inclusion of the $L_2$ regularization term in~\eqref{coupled_reg} can be interpreted as a prior expectation that the norm of the inferred operators should not be arbitrarily large. We focus on this second viewpoint and establish the connection between the $L_2$ regularization term in either LS problem of~\eqref{coupled_reg} and the loci of the inferred linear operator eigenvalues via the Gershgorin disk theorem~\cite{gershgorin}.

For ease of exposition, we consider an autonomous linear system as in~\eqref{fom_linear}. We first briefly state the Gershgorin disk theorem for the linear operator $\mathbf{A}$. Denoting the entries of a matrix $\mathbf{A} \in \mathbf{R}^{n \times n}$ as $\mathbf{A}_{i,j}$, an upper bound on the locus of the $k$-th eigenvalue of $\mathbf{A}$ is given by 

\begin{equation}
	\label{gersh_theorem}
	\left|\lambda_k - \mathbf{A}_{i,i}\right| \leq \sum_{j \neq i} \left|\mathbf{A}_{i,j} \right|, \qquad \text{for some} \; i \in \{1, \dots, n \}.
\end{equation}

The eigenvalue $\lambda_k$ lies in at least one of the Gershgorin disks, with the disk center being the diagonal entry $\mathbf{A}_{i,:}$ and the disk radius being the r.h.s term of~\eqref{gersh_theorem}. By relaxing the upper bound of~\eqref{gersh_theorem} we get  

\begin{equation}
	\label{gersh_2}
	\left|\lambda_k - \mathbf{A}_{i,i}\right|  \leq \sum_{j \neq i} \left|\mathbf{A}_{i,j} \right| \leq \sum_j \left|\mathbf{A}_{i,j} \right| = \underbrace{\left\Vert \mathbf{A}_{i,:} \right\Vert_1 \leq \sqrt{n} \left\Vert \mathbf{A}_{i,:} \right\Vert_2}_{\text{Cauchy-Schwarz inequality}}, \qquad \text{for some} \; i \in \{1, \dots, n \}.
\end{equation}

\noindent From an inference viewpoint,~\eqref{gersh_2} shows that by penalizing the $L_2$ norm of the rows of $\mathbf{A}$ in~\eqref{coupled_reg}, we restrain an upper bound of the corresponding Gershgorin disk radii where the eigenvalues of $\mathbf{A}$ lie. We write \eqref{gersh_2} for each linear system inference task in \eqref{coupled_inference}. For OpInf, we get

\begin{equation}
	\label{opinf_ineq}
	\left|\lambda_k - \mathbf{A}_{{\text{RR}}_{i,i}}\right|   \leq \sqrt{r} \left\Vert \mathbf{A}_{{\text{RR}}_{i,:}} \right\Vert, \qquad \text{for some} \; i \in \{1, \dots, r \},
\end{equation}

\noindent while for sFOM, \eqref{gersh_2} is written as

\begin{equation}
	\label{sfom_ineq}
	\left|\lambda_k - \mathbf{A}_{{\text{FF}}_{i,i}}\right|   \leq \sqrt{n_{Q_i}} \left\Vert \mathbf{A}_{{\text{FF}}_{i,Q_i}} \right\Vert, \qquad \text{for some} \; i \in \{1, \dots, n_F \}.
\end{equation}

This establishes a connection between the Gershgorin disks radii of the inferred OpInf and sFOM linear operator in \eqref{opinf_ineq} and \eqref{sfom_ineq} and the $L_2$ regularization terms employed for performing inference via~\eqref{coupled_reg}. Based on this observation, we propose a variant of the $L_2$ regularization that promotes the stability of the inferred OpInf and sFOM linear operator and has a closed-form solution.

\subsubsection{Gershgorin regularization}
\label{subs:gershreg}

A linear dynamical system, as~\eqref{fom_linear}, is asymptotically stable if all the eigenvalues of $\mathbf{A}$ lie in the left semi-plane, i.e.,

\begin{equation}
\label{eigs_condition}
	\text{Re}(\lambda_k)<0, \qquad \forall \; k \in \{1, \dots, n \}.
\end{equation}

Although the $L_2$ regularization~\eqref{coupled_reg} has been successful in many applications, inferred dynamical systems are often unstable. This can be explained in terms of~\eqref{gersh_2} for an autonomous linear system~\eqref{fom_linear}. The $L_2$ regularization of~\eqref{coupled_reg} restrains the loci of eigenvalues $\lambda$ towards the origin of the complex plane, with no directionality towards the left semi-plane. For this reason, many studies have proposed modifications of the LS problem to obtain stable non-intrusive linear ROMs or FOMs by regularization scaling and post-processing~\cite{SAWANT2023115836}, constrained optimization~\cite{prakash2024datadriven} or a nonlinear LS formulation~\cite{goyal2023guaranteed}. 

We leverage expression~\eqref{gersh_2} to propose a regularization strategy for each inference task in~\eqref{coupled_inference}, that promotes stability of the inferred linear operator, while admitting a closed-form solution. Essentially, we would like not only to restrain the Gershgorin disk radii but also to drive the disk centers to the left complex semi-plane. In this way, we promote eigenvalues with a negative real part, as dictated by \eqref{eigs_condition}. It should be noted that the strong stability enforcement by the Gershgorin theorem as in \cite{prakash2024datadriven} is a sufficient but not necessary stability condition for a linear operator $\mathbf{A}$.

We follow the notation used in sFOM~\eqref{sfom_ls}, and note that the exposition for OpInf in~\eqref{opinf} follows analogously. Specifically, we focus on~\eqref{coupled_structure_linear}, and complement the $L_2$ regularization for the sFOM problem with one more term which penalizes the diagonal entries of $\mathbf{A}$. By regularizing the $L_2$ norm of the LS solution, we penalize the norm of the $i$-th row of the linear operator, which accounts for restricting the upper bound of the $i$-th Gershgorin disk radius, according to~\eqref{sfom_ineq}. Adding a penalization term for the inferred diagonal entry $\mathbf{A}_{i,i}$ accounts for driving the Gershgorin disk center towards the left complex semi-plane according to~\eqref{gersh_2}, thus introducing a ``directionality" of the regularization for the loci of the linear operator eigenvalues towards the left semi-plane. The resulting LS formulation is written as

\begin{equation}
	\label{sfom_gersh}
	\min _{\bm{\beta}_i}{\left( \left\|\bm{\beta}_i^\top {\bm{\mathcal{D}}_{\text{F}_i}} -  \left. \frac{ \textrm{d}\mathbf{X}_\text{F}}{\textrm{dt}} \right\vert_{i} \right\|^2_{2} + \eta_{1_i} \left\|\bm{\beta}_i \right\|^2_2 + \eta_{2_i} {\bm{\beta}_i}_m \right)}, \qquad
	\bm{\beta}_i = 
	\begin{bmatrix}
		\mathbf{A}_{\text{FF}_{i, Q_i}}^\top \\ \mathbf{A}_{\text{FR}_{i, :}}^\top
	\end{bmatrix},
\end{equation}

\noindent where ${\bm{\beta}_i}_m =\mathbf{A}_{\text{FF}_{i,i}}$ and $\eta_{2_i}$ is a corresponding weighting factor. 

Formulation~\eqref{sfom_gersh} can be easily incorporated into existing LS solvers since it admits the closed-form solution 

\begin{equation}
	\label{soln_gersh}
	\left({\bm{\mathcal{D}}_{\text{F}_i}} \bm{\mathcal{D}}_{\text{F}_i}^\top + {\eta_{1_i}} \mathbb{I} \right) {\bm{\beta}_i} =  {\bm{\mathcal{D}}_{\text{F}_i}} \left. \frac{ \textrm{d}\mathbf{X}_\text{F}}{\textrm{dt}} \right\vert_{i}^\top - \eta_{2_i} \mathbf{e}_m,
\end{equation}

\noindent where $\mathbf{e}_m$ is the unitary vector such that $\mathbf{e}_m^T {\bm{\beta}_i} = \mathbf{A}_{\text{FF}_{i,i}}$. 

\begin{figure}[!htb]
	\centering
	\includegraphics[width=.8\textwidth]{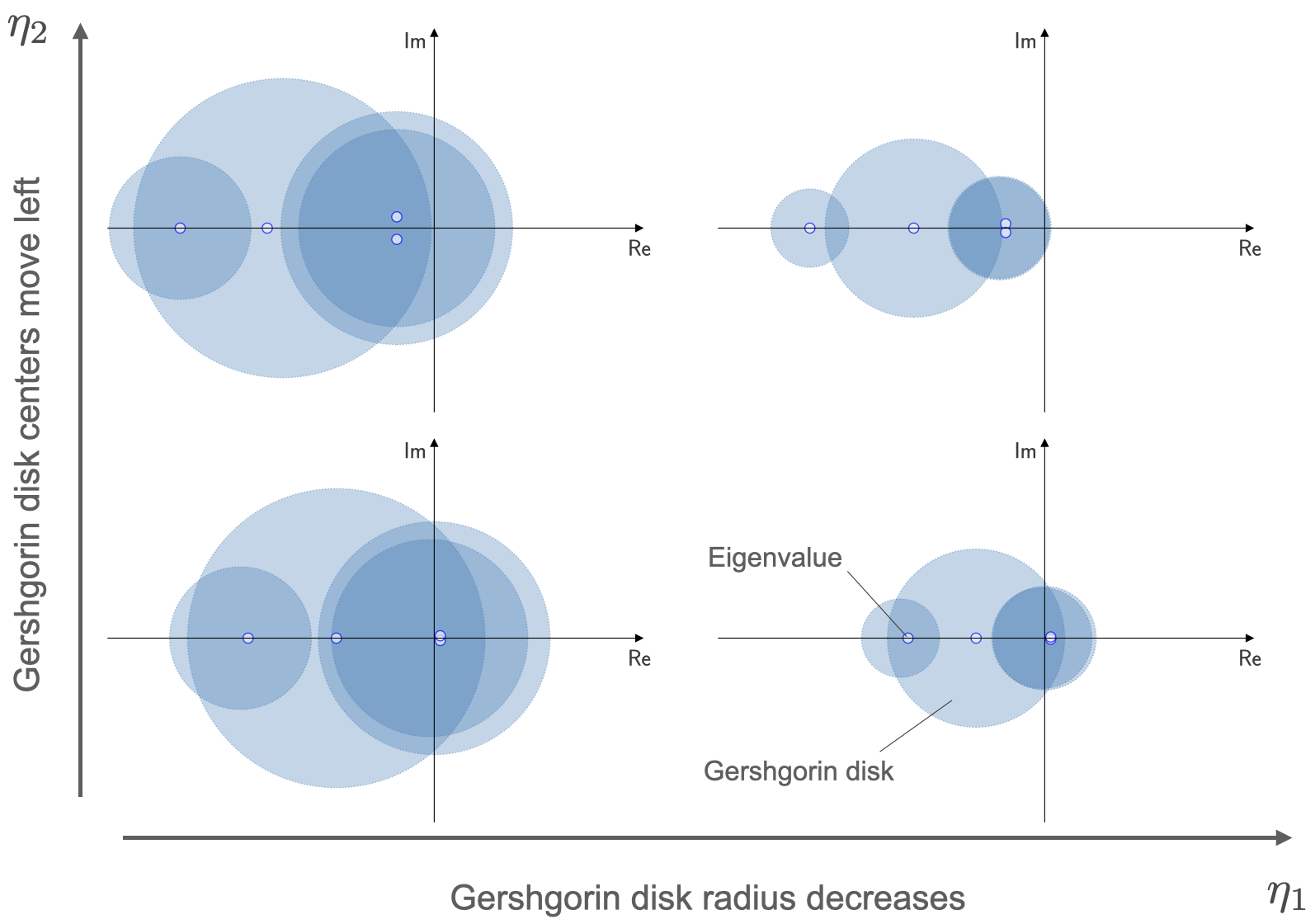}
	%\includegraphics[width=0.32\textwidth, trim=13 0 13 0, clip]{figures/no_reg.pdf}
	%\hfill
	%\includegraphics[width=0.32\textwidth, trim=13 0 13 0, clip]{figures/10_0_no_gersh_reg.pdf}
	%\hfill
	%\includegraphics[width=0.32\textwidth, trim=13 0 13 0, clip]{figures/gersh_regularization.pdf}
	\caption{The effect of the Gershgorin $L_2$ regularization on the eigenvalues of an inferred linear operator. From left to right, we increase the factor $\eta_1 \geq 0$, corresponding to Tikhonov (or $L_2$) regularization. This can be interpreted as a stronger penalization of the Gershgorin disk radii of the inferred linear operator. From bottom to top, we increase the factor $\eta_2 \geq 0$ and promote negative values for the diagonal entries of the linear operators. The introduced Gershgorin regularization in \eqref{soln_gersh} combines both $\eta_1$ and $\eta_2$ terms.}
 \label{fig:gersh}
\end{figure}

Considering $\eta_{2_i}=0$ results in the typical Tikhonov regularization in~\eqref{coupled_reg}. By introducing $\eta_{2_i} > 0$, formulation~\eqref{sfom_gersh} is a stability promoting $L_2$ regularization variant which encourages the inference of linear operators with eigenvalues on the left complex semi-plane. A schematic representation of the effect of terms $\eta_1$ and $\eta_2$ to the inferred linear operator eigenvalues is given in \Cref{fig:gersh}. 

For nonlinear problems, the expositions in~\eqref{sfom_gersh} and~\eqref{soln_gersh} hold for the system's linear operator (the entries of $\mathbf{e}_m$ are $0$ for the nonlinear terms). For example, in quadratic systems, the stability of the linear term is promoted via~\eqref{sfom_gersh} and the quadratic term is regularized as in~\cite{Gkimisis2023}, since regularizing the norm of the quadratic matrix is known to promote stability~\cite{Kramer2021}. As indicated by~\Cref{fig:gersh}, depending on the $\eta_1$ and $\eta_2$ values, the proposed regularization encourages stability of the linear operator, thus allowing to also accommodate dynamical systems with inherently unstable linear operators~\cite{Noack2003}. Furthermore, we note that stability of both the sFOM and the OpInf models in~\eqref{coupled_structure_linear} promotes, but does not guarantee, stability of the two-way coupled system (see \Cref{fig:matrix_structures}). However, in practice, Gershgorin regularization has proven to be valuable for assisting the inference of both robust sFOM, as well as OpInf models. This is showcased in both numerical test cases in \Cref{sec:results}.

%seperate section, computational cost analysis
\section{Computational cost analysis}
\label{sec:cost}
Since the OpInf-sFOM approach includes full-order information to build a locally accurate ROM, we investigate its offline and online computational cost. Although spatially localized, the sFOM inference is performed at the full-order level. Hence, we present a quantitative parametric analysis of the offline and online computational efficiency of the coupled OpInf-sFOM method with respect to the size of the sFOM subdomain and the achievable reduction through OpInf. This analysis reveals the applicability and efficiency of the presented approach for the inference of a given system.

\subsection{Asymptotic offline computational cost}
\label{offline_cost}

% For both OpInf and sFOM, the inference task, i.e., solving problems~\eqref{opinf} and~\eqref{sfom_ls}, is the main contributor to the inference offline computational cost. 

The main contributor for the sFOM-OpInf offline cost is solving the OpInf and sFOM LS inference problems in \eqref{coupled_inference} which asymptotically scale as $\mathcal{O}(n_T m^2)$, where $n_T$ is the number of training snapshots and $m$ is the number of unknowns. This assumes that there are more equations than unknowns $(n_T \gg m)$ ~\cite{strang2021introduction}. 
For the sFOM inference in~\eqref{sfom_ls}, we consider $s$ adjacent DOFs for any DOF $i$ and a model polynomial order $k$. Then the number of unknowns $m$ scales as $m \propto \binom{s+1}{k}$, with a leading polynomial term equal to $\frac{s^k}{k !}$. Based on the above, we conclude that the leading term for the sFOM offline cost is 

\begin{equation}
\label{sfom_cost_1}
C^{\text{off}}_{\text{sFOM}} \approx n_\text{F}\left( \frac{s^{k}}{{k !}} \right)^2 n_T,
\end{equation}

\noindent since $n_\text{F}$ LS problems need to be solved in~\eqref{sfom_ls}. 

For OpInf, we consider the LS problem in~\eqref{opinf} where $m$ scales as $m \propto \binom{r+n_\text{I}+1}{k}$. Analogously to the sFOM computation, the OpInf offline computational cost scales as

\begin{equation}
\label{op_cost_2}
C^{\text{off}}_{\text{OpInf}} \approx r \left( \frac{\left( r+ n_\text{I} \right) ^k}{{k !}} \right)^2 n_T,
\end{equation}

\noindent since the information from the $n_\text{I}$ sFOM DOFs are encoded as an input term to the OpInf model. The corresponding cost of a global OpInf model with dimension $r_g$ and the cost of a global sFOM for all $n$ DOFs would scale as

%Comparing~\eqref{op_cost_1} and~\eqref{op_cost_2}, we expect that the offline computational cost in~\eqref{op_cost_2} under the implicit formulation in~\eqref{implicit_opinf} will be higher due to the treatment of the coupling term. We will focus on this formulation as a worst-case scenario. 

\begin{equation}
\label{op_glob_1}
C^{\text{off}}_{\text{OpInf}_g} \approx \frac{r_g^{2k+1}}{{k !}^2} n_T, \qquad C^{\text{off}}_{\text{sFOM}_g} \approx n \left( \frac{s^{k}}{{k !}} \right)^2 n_T,
\end{equation}

To compare how the computational cost of the coupled OpInf-sFOM approach scales with respect to a global sFOM or a global OpInf model, we take the ratio of the corresponding expressions. This results in

\begin{equation}
\label{ratio_opg}
\frac{C^{\text{off}}_{\text{OpInf}-\text{sFOM}}}{C^{\text{off}}_{{\text{OpInf}}_g}} \approx \left [ \frac{n_{\text{F}}}{r} \left( \frac{s}{r}\right)^{2k}  +  \left(1+\frac{n_\text{I}}{r} \right)^{2k} \right ] \left(\frac{r_g}{r}\right)^{-2k-1},
\end{equation}

\noindent and correspondingly 

\begin{equation}
\label{ratio_sfomg}
\frac{C^{\text{off}}_{\text{OpInf}-\text{sFOM}}}{C^{\text{off}}_{{\text{sFOM}}_g}} \approx \frac{n_{\text{F}}}{n} + \frac{r}{n}\left(1+\frac{n_\text{I}}{r} \right)^{2k} \left(\frac{s}{r}\right)^{-2k}.
\end{equation}

The $n_\text{I}$ DOFs close to the OpInf-sFOM interface scale as $n_\text{I} \in \mathcal{O} \left( n_{\text{F}}^{d-1} \right)$, where $d$ is the spatial dimension of the problem (e.g., $d=3$ for a 3D problem). We can then compute the offline computational cost ratios in~\eqref{ratio_opg} and~\eqref{ratio_sfomg} for families of problems with a given polynomial nonlinear order $k$ and a given spatial dimension $d$. We consider typical values for $s$ and $r$, with $s \in \mathcal{O} \left(3^d \right)$ (for a 3-point stencil in 1D) and $ r \in \mathcal{O}(10)$.

\begin{figure}[!htb]
\centering
\includegraphics[width=0.48\textwidth, trim=0 0 35 35, clip]{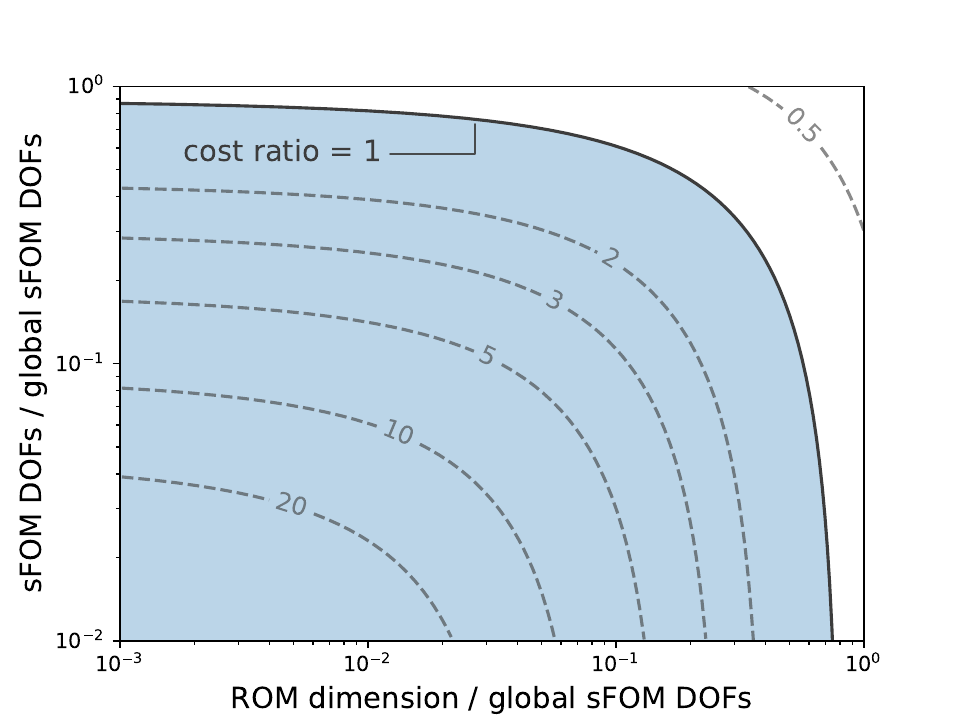}
\hfill
\includegraphics[width=0.48\textwidth, trim=0 0 35 35, clip]{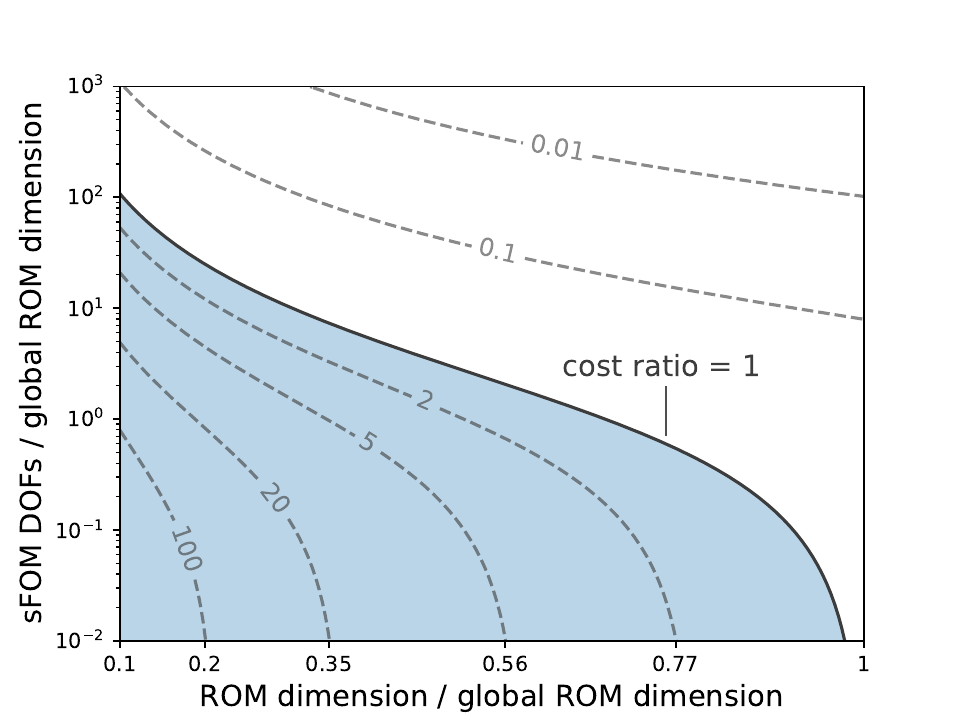}
\caption{Asymptotic offline speedup estimates for linear 2D problems. Left panel: contours of offline cost ratio of global sFOM over coupled OpInf-sFOM approach. Right panel: contours of offline cost ratio of global OpInf over coupled OpInf-sFOM approach. The coupled OpInf-sFOM cost decreases for a contracting FOM subdomain and a reducible ROM subdomain.}
\label{fig:cost}
\end{figure}

In \Cref{fig:cost}, the asymptotic offline speedup of the coupled approach in~\eqref{coupled_inference} with respect to the global sFOM cost in~\eqref{sfom_ls} and the global OpInf cost in~\eqref{opinf} is given for linear 2D problems. The main observation in both graphs is the increasing cost benefit of the coupled approach for smaller FOM subdomains and smaller ROMs. The coupled approach is more cost-effective if the dynamics related to slow singular value decay are localized while the dynamics on the rest of the domain can be embedded into a low-dimensional space. \eqref{ratio_opg} and \eqref{ratio_sfomg} provide insight into the cost benefit of the coupled approach for any given problem parameters ($k, d, s, r$). Although this analysis does not consider the prediction accuracy, the quantities in \Cref{fig:cost} can be estimated by problem-specific criteria. For example, a value for the dimension ratio of global to local OpInf basis -- on the $x$-axis of the right panel in \Cref{fig:cost} -- could be estimated from the corresponding, snapshot data normalized singular value decay, for a given projection error threshold. 

\subsection{Online cost considerations}
\label{online_cost_pred}

During the online phase, the computational cost is due to the matrix-vector multiplications for the computation of the time derivative in~\eqref{coupled_structure_quad} and the selected time integration scheme. To estimate the online cost speedup of the OpInf-sFOM in comparison to a FOM code, we analyze the case in which the model integrated in time by the FOM code is explicitly known and of the form in~\eqref{fom_structure}, and the operators have the same sparsity pattern as that employed for OpInf-sFOM. In practice, for many demanding applications (e.g., in CFD codes) the operators are not explicitly formed, while features such as limiters, look-up tables, and iterative and stabilization schemes are employed~\cite{farcas2023improving}. We also consider the case where the FOM code uses the same time integration scheme as the OpInf-sFOM model in~\eqref{coupled_structure_quad}. Though the implementation of the FOM code may practically vary from the setup analyzed here, the above considerations provide a reasonable baseline for an online computational cost comparison.

Following the same asymptotic analysis as in \Cref{offline_cost} for the leading term of the matrix-vector multiplication cost, we get an estimate of the scaling of the online sFOM-OpInf cost as well as the FOM cost as

\begin{equation}
\label{online_sFOM_OpInf}
C^{\text{on}}_{\text{sFOM}} \approx n_\text{F} \frac{s^{k}}{{k !}} n_t, \qquad
C^{\text{on}}_{\text{OpInf}} \approx r \frac{\left( r+ n_\text{I} \right) ^k}{{k !}}  n_t, \qquad C^{\text{on}}_{\text{FOM}} \approx n \frac{s^{k}}{{k !}} n_t.
\end{equation}

 We can then write the ratio of the OpInf-sFOM over the FOM online cost as

\begin{equation}
\label{online_cost_ratio}
\frac{C^{\text{on}}_{\text{OpInf-sFOM}}}{C^{\text{on}}_{\text{FOM}}} \approx \frac{r}{n} \left( \frac{r+n_\text{I}}{s}\right)^k + \frac{n_\text{F}}{n}.
\end{equation}

\begin{figure}[!htb]
\centering
\includegraphics[width=0.8\textwidth]{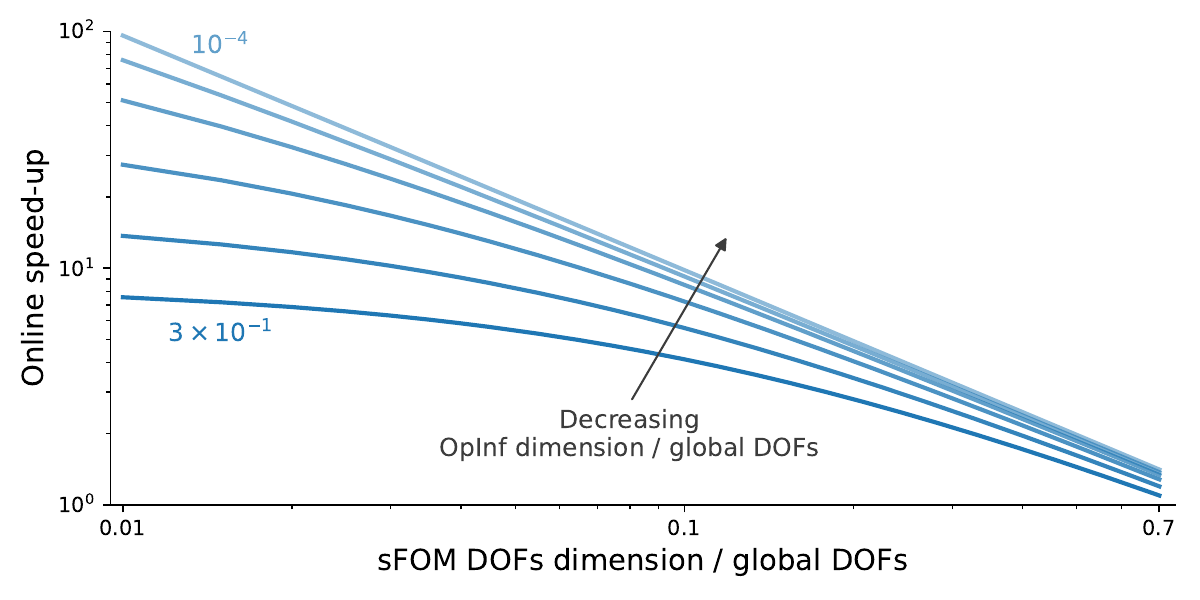}
\caption{OpInf-sFOM online speedup compared to FOM, with respect to sFOM subdomain size and OpInf dimension, for 2D, linear problems with a $3$-point numerical stencil along each spatial coordinate. As the sFOM subdomain size increases, the computational gain decreases toward 1. For smaller OpInf dimensions, the obtained speedup is increased.}
\label{fig:online}
\end{figure}

The expression in~\eqref{online_cost_ratio} indicates an almost inversely proportional relation between the achieved online speedup and the relative size of the sFOM subdomain $n_\text{F}/n$. This is intuitively expected since computations on the sFOM subdomain are performed at the FOM level. An additional term arises from the OpInf model, where coupling sFOM DOFs also appear. For given problem parameters $d, \; k$ and selecting typical values for $s$ and $r$, we can draw curves on the expected speedup with respect to the sFOM size and the non-dimensional reduction achieved for OpInf ($r/n$). \Cref{fig:online} shows the curves for linear ($k=1$) two-dimensional ($d=3$) problems.

\section{Numerical results}
\label{sec:results}

We apply the coupled OpInf-sFOM formulation in~\eqref{coupled_inference} using the Gershgorin regularization introduced in~\Cref{subs:gershreg} for two test cases. The first is an illustrative example of a 1D Burgers' equation with periodic boundary conditions, for which a coupled, quadratic OpInf-sFOM model is inferred. The second test case is a weakly nonlinear 2D model for the ice thickness dynamics of the Pine Island Glacier.

\subsection{1D Burgers' equation}
\label{subs:1d_burg}

The 1D Burgers' equation with periodic boundary conditions reads as
\begin{equation}
\left \{ 
	\begin{aligned}
		\label{1Dburg}
		\frac{\partial w}{\partial t} - c w \frac{\partial w}{\partial z}  -\nu \frac{\partial^2 w}{\partial z^2} &=0,\\
		w(z=0) & = w(z=L), \\
		w(t=0) & = w_0(z). \\
	\end{aligned}
	\right.
\end{equation}

\noindent where $c$ denotes the advection speed, $\nu$ denotes the kinematic viscosity and $L$ is the length of the domain along $z$. We discretize~\eqref{1Dburg} in space and time with $\Delta z=2 \times 10^{-2}$, $\Delta t=2.5 \times 10^{-2} s$ and choose $c=0.5$, $\nu = 10^{-2}$, $L=10$, and overall simulation time of $T=$ \SI{18}{s}. The initial condition is $w_0(z) = \exp \left\{-\frac{{(z-L/2)}^2}{1.2} \right\}+0.1 \cos(2\pi z/L)+0.1 \cos(4\pi z/L)$, corresponding to a Gaussian centered at $z=L/2$ with periodic small-amplitude disturbances and therefore the main dynamical feature of the system is a traveling wave towards $z=L$. The simulation is performed via the Fast Fourier Transform (FFT) method in~\cite{git1dburg} and a data matrix $\mathbf{X} \in \mathbb{R}^{500 \times 720}$ is assembled.

Due to the selected initial condition, the region $z \in \left[ 5, 10 \right]$ is expected to exhibit transport-dominated dynamics, and thus a slow singular value decay. Indeed, using the training data for the first \SI{9}{s}, we observe a significant gap in the singular value decay rate between the data in $z \in \left[ 0, 5 \right)$ and those in $z \in \left[ 5, 10 \right]$. We therefore accordingly assign these subdomains to OpInf and sFOM. We proceed to apply the OpInf-sFOM~\eqref{coupled_inference} with the Gershgorin regularization in~\eqref{sfom_gersh}. The underlying discretized model for~\eqref{1Dburg} has a quadratic nonlinearity due to the advection term (equation \eqref{fom_structure} with $\mathbf{T}$ and all higher order terms equal to zero). Hence, the coupled OpInf-sFOM structure is the one in~\eqref{coupled_structure_quad}, with quadratic nonlinearity and $\mathbf{B}=\mathbf{0}$, $\mathbf{c}=\mathbf{0}$.

\begin{figure}[!htb]
    \centering
    \begin{minipage}{0.48\textwidth}
        \centering
        \includegraphics[width=0.99\textwidth]{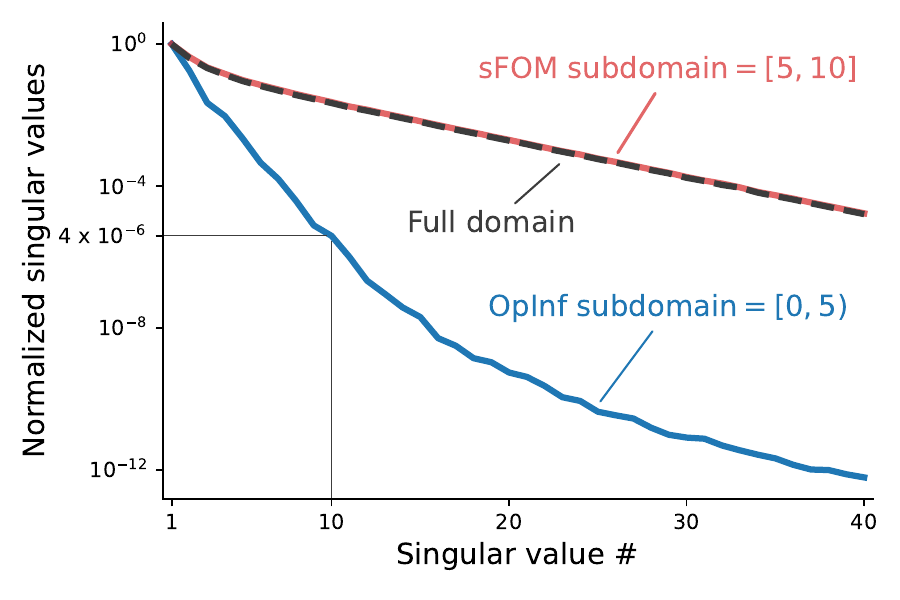}
\caption{Normalized singular value decay over training data. The training data indicate a significant difference in the singular value decay rate between the OpInf subdomain $x \in \left[ 0, 5 \right)$ and the sFOM subdomain $x \in \left[ 5, 10 \right]$.} 
\label{fig:svals_burg}
    \end{minipage}\hfill
    \begin{minipage}{0.48\textwidth}
        \centering
        \includegraphics[width=0.99\textwidth]{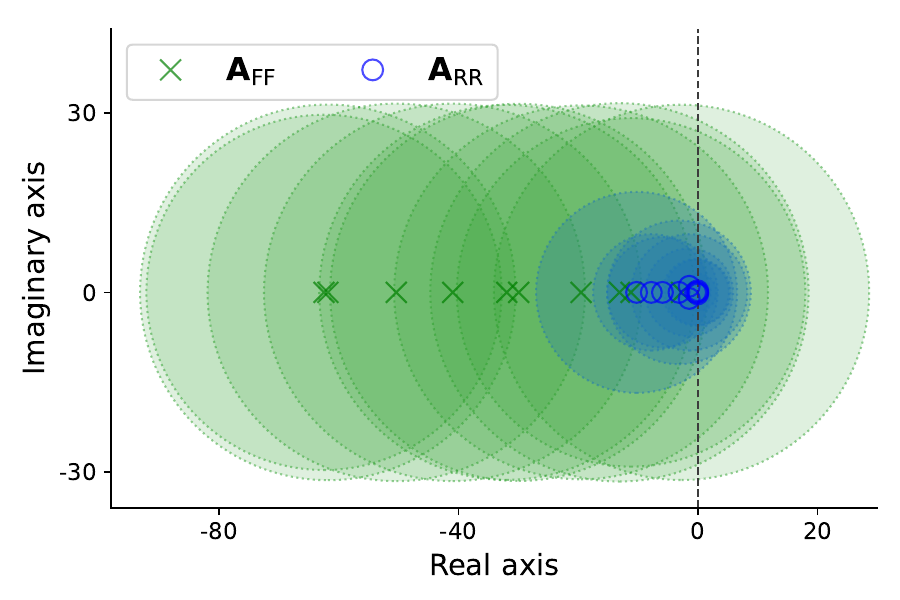}
\caption{Eigenvalues and Gershgorin disks for the inferred OpInf and sFOM linear operators $\mathbf{A}_{\text{FF}}$ and $\mathbf{A}_{\text{RR}}$ in \eqref{coupled_inference}. Employing the introduced regularization in \eqref{sfom_gersh} for the linear OpInf and sFOM terms} leads to stable linear operators. We observe that the sFOM inference retrieves a symmetric negative definite matrix $\mathbf{A}_{\text{FF}}$, which is in line with the underlying, diffusive linear operator (see \eqref{1Dburg}).
\label{fig:eigs_burg}
    \end{minipage}
\end{figure}

For the OpInf subdomain in $z \in \left[ 0, 5 \right)$, a reduced dimension $r=10$ captures more than $99.9\%$ of the system energy and is sufficient for an accurate ROM. We compute the optimal regularization values for OpInf via an L-curve criterion, by testing $20$, logarithmically spaced values $\eta_1^{\text{OpInf}} \in \left[ 10^{-3}, \; 10^{0}\right]$. Since the linear and quadratic operators scale differently~\cite{mcquarrie2021data}, we penalize the quadratic term with a factor of $200 \times \eta_1^{\text{OpInf}}$. For the Gershgorin penalization in~\eqref{sfom_gersh}, we find that a factor of $\eta_2^{\text{OpInf}} = 0.05 \times \eta_1^{\text{OpInf}}$ produces accurate results, and thus we compute only one L-curve with respect to the values of $\eta_1$. All eigenvalues of the resulting $\mathbf{A}_{\text{RR}}$ are on the left complex semi-plane, as indicated by \Cref{fig:eigs_burg}.

\begin{figure}[!htb]
\centering
%\includegraphics[width=0.48\textwidth]{figures/1d_burgers_svd.pdf}
%\hfill
\includegraphics[width=1.\textwidth, trim=80 0 50 0, clip]{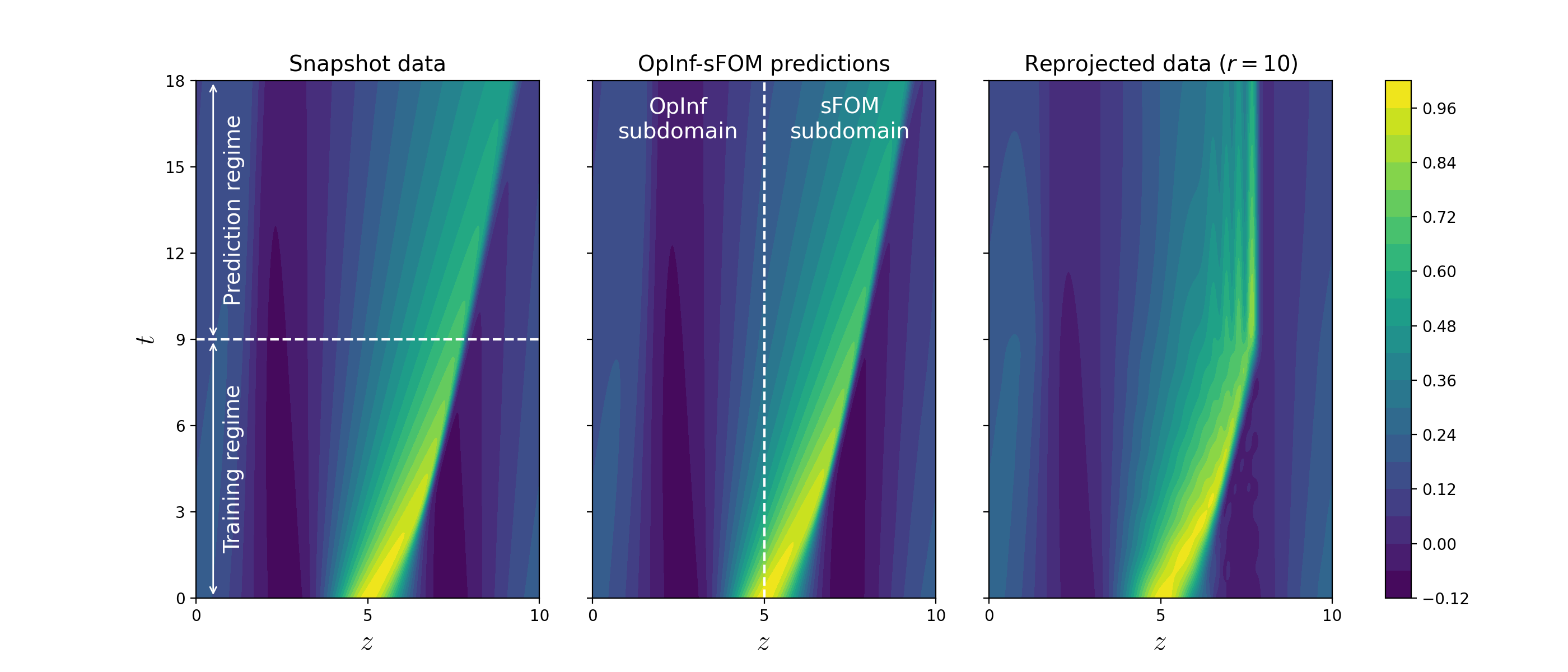}
\caption{OpInf-sFOM predictions for a 1D Burgers' equation test case. Left panel: simulation data obtained via the FFT method. We use the first \SI{9}{s} as training data. Central panel: state predictions via OpInf-sFOM. Accurate predictions for the system dynamics are made beyond the span of the training snapshot matrix. Right panel: projected data to global basis $\mathbf{\Phi} \in \mathbf{R}^{1000 \times 10}$ constructed by the truncated SVD of the training data ($t \leq 9 s$) which retains more than $99.9\%$ of the system energy. A classic projection-based approach is bound to the span of $\mathbf{\Phi}$ and cannot predict the discontinuity propagation.}
\label{fig:1d_burg}
\end{figure}

For the sFOM subdomain in $z \in \left[ 5, 10 \right]$, we use a sparsity pattern with a $3$-point stencil. Since the 1D mesh is uniform, we solve~\eqref{sfom_ls} for each $i$ by concatenating data from $5$ randomly selected sFOM DOFs. Analogously to OpInf, we choose an optimal regularization parameter via an L-curve criterion, by testing $20$ logarithmically spaced values $\eta_1^{\text{sFOM}} \in \left[ 10^{-3}, \; 10^{-8}\right]$. A penalization of the quadratic term with values $10 \times \eta_1^{\text{sFOM}}$ and a penalization of the linear operator diagonal entries with values $\eta_2^{\text{sFOM}} = 50 \times \eta_1^{\text{sFOM}}$ are set. After the optimization of $\eta_1$ via the L-curve, the resulting linear operator $\mathbf{A}_{\text{FF}}$ is stable. The eigenvalues of $\mathbf{A}_{\text{FF}}$ and the corresponding Gershgorin disks in~\eqref{gersh_theorem} are given in \Cref{fig:1d_burg}.

The predictions made by the OpInf-sFOM model are presented in \Cref{fig:1d_burg}. Even though the discontinuity forms and propagates in a region beyond the one of the training regime, inferring an sFOM for $z \in \left[ 5, 10 \right]$ allows for an accurate prediction of the dynamics. On the contrary, projecting the FOM to a global basis computed by the training regime data is by definition insufficient for capturing this advective pattern, since it lies beyond the span of the training snapshot matrix $\mathbf{X}$. In parallel, inferring an OpInf model for $z \in \left[ 0, 5 \right)$ allows for a $25\times$ reduction of the DOFs on that subdomain. Therefore, the online speedup factor is $\sim 1.22$ in \eqref{online_cost_ratio} relative to the FOM code. A parametric study of the effect of the OpInf-sFOM interface position on the properties of the inferred, coupled system is given in \ref{appendix:app}.

\subsection{Ice thickness dynamics}
In our second test case, we employ the proposed OpInf-sFOM formulation to model the Pine Island Glacier ice thickness dynamics, and evaluate the method's predictive capabilities for large scale applications.
As the fastest melting glacier in Western Antarctica \cite{rignot2013ice}, the Pine Island Glacier is of special interest for ice melt predictions and the estimation of future sea level rise.
However, numerical predictions are challenging due to the involved parametric uncertainties, the coupling of different physical dynamics (such as ice thickness and velocity), the fine resolution required for resolving them, and the consequent high computational cost \cite{greve2009dynamics}.
In the following, we use our OpInf-sFOM formulation to construct a ROM of the ice thickness $h$ to enable fast predictions of future ice mass loss.
The ROM is parameterized by the ice melting rate $\alpha$ at the ice-ocean interface to account for uncertain environmental circumstances. 
In particular, the sFOM part of the ROM resolves the expanding floating portion of the ice where the oceanic melt is applied.

% To mitigate these issues and to allow multi-fidelity uncertainty quantification (UQ), we construct a parameterized sFOM-OpInf ROM for the ice thickness $h$ of the PIG for varying ice melting rates $\alpha$ at the ice-ocean interface.

\subsubsection{Constitutive equations}
We model the ice thickness $h(t)$, $t \in [0, 20]$, for 20 years on the domain $\Omega \subset \mathbb{R}^2$ of the Pine Island Glacier through
\begin{equation}
\label{diffusion}
    \frac{ \textrm{d}h}{\textrm{dt}} = -\nabla \cdot h \mathbf{v} + \nabla \cdot (D \nabla h) + m_s - m_b(h, \alpha), \qquad \text{in } \Omega, 
\end{equation}
with homogeneous Neumann boundary conditions $
\nabla h \cdot \mathbf{n} = 0$ on $\partial \Omega$.
In~\eqref{diffusion}, $D \in \mathbb{R}^{2\times 2}$ is a symmetric positive definite diffusivity term introduced for numerical stabilization, $m_s$ [\si{m/yr}] describes the rate of mass deposition on the ice surface, $m_b(h, \alpha)$ [\si{m/yr}] describes melting at the ice base, and $\mathbf{v}$ [\si{m/yr}] is a depth-averaged ice velocity field.

For this test case, $\mathbf{v}$ and $m_s$ are constant in time. From a physical point of view, this assumption is reasonable due to the observed slow ice velocity changes and the usage of a time-averaged surface mass balance. Instead, we focus on the nonlinearity introduced by the local basal melt forcing term $m_b(h, \alpha)$:
In our model, basal melt occurs at the interface with the (warmer) ocean water.
Since the ice density $\rho_{\rm{ice}}$ is lower than the sea water density $\rho_{\rm{water}}$, the ice shelf floats at the water surface when $h$ is thinner than the floating height $h_f$, i.e., when
\begin{equation}\label{eq:floating-condition}
    h(z) < h_f(z) := -\frac{\rho_{\text{water}}}{\rho_{\text{ice}}} b(z) 
\end{equation}
with ocean water density $\rho_{\text{water}} = \SI{1023}{kg/m^3}$, ice density $\rho_{\text{ice}} = \SI{917}{kg/m^3}$, and bedrock elevation field $b$ from \cite{nitsche2007bathymetry}.
At these positions, $m_b(h, \alpha)$ applies the ice melt rate $\alpha$, i.e.,
\begin{equation}
\label{basal_melt}
m_b(h, \alpha) = \alpha \; \mathds{1}_{\{h(z) > h_f(z)\}},
% m_b(h, \alpha) = \alpha \delta \left(h(\mathbf{x})>h_f(\mathbf{x})\right), 
\end{equation}
When~\eqref{eq:floating-condition} does not hold, the ice is pressed against the (cold) bedrock and no melting occurs.
As ice becomes thinner, e.g., due to ice melting or sliding, previously grounded areas of the ice shelf detach and start to float and consequently melt.
The accurate estimation of the floating area is of major interest in predictive simulations as it determines the ice shelf's physical behavior. 
We then expect, by the aforementioned physical mechanisms, that the region where the nonlinear forcing $m_b(h, \alpha)$ is active will also exhibit transport-dominated dynamics.

% \noindent with $h_t(\mathbf{x})$ being a space-dependent threshold value for the ice thickness. Equation~\eqref{basal_melt} introduces a localized nonlinearity to~\eqref{diffusion}.

\begin{figure}[!t]
\centering
  \includegraphics[width=.65\textwidth]{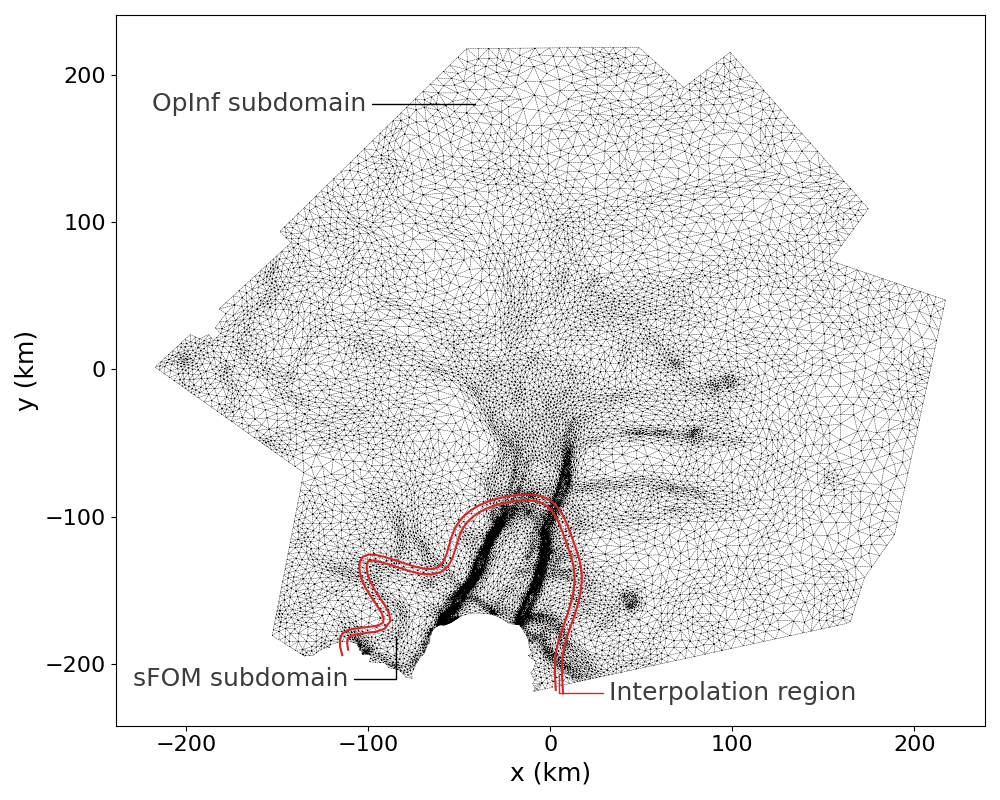}
\caption{Pine Island Glacier mesh. We collect data from a mesh with $12,494$ DOFs. For the inference, we decompose the domain with $4,585$ DOFs on the sFOM subdomain and $7,909$ DOFs on the OpInf subdomain. As a post-processing step, the solution is interpolated along a narrow strip along the ROM-FOM interface.}
\label{fig:ice1}
\end{figure}

We discretize~\eqref{diffusion} in space using linear finite elements ($12494$ DOFs, see Figure \ref{fig:ice1}).
The mesh was generated from an adapted outline of the Pine Island Glacier in \cite{seroussi2014sensitivity} to approximate the ice velocity in \cite{rignot2011ice, rignot2011measures} at a resolution between \SI{250}{m} and \SI{10}{km}.
This resolution was chosen according to study \cite{seroussi2014sensitivity}.
For the surface accumulation $m_s$, we use \cite{vaughan1999reassessment} as a yearly average.
To obtain the ice velocity $\mathbf{v}$, we solve the Shallow Shelf Approximation (\cite{greve2009dynamics}, Chapter 6.3) with basal friction field inferred from the surface velocity measurements in \cite{rignot2011ice, rignot2011measures}.
For the initial condition $h(0)$, we subtract the surface topography \cite{bamber2009new, griggs2009newPart2} from the ice bed \cite{nitsche2007bathymetry} such that $t=0$ roughly corresponds to the year 2004. 
Given a melt rate $\alpha$, the model~\eqref{diffusion} is solved until year $t=20$ using the implicit Euler method with step size $\Delta t = 0.01$. 
On a single core, a solve takes approximately \SI{9}{min}.
The model is implemented in the Ice Sheet System Model (ISSM) \cite{larour2012continental}.

\subsubsection{OpInf-sFOM predictions}

We exploit the physical insight from equations~\eqref{diffusion} and~\eqref{basal_melt} for the coupled OpInf-sFOM model. Knowing that~\eqref{diffusion} is linear and that the nonlinear input term in~\eqref{basal_melt} is localized, we infer a linear OpInf model with $r=10$ retaining more than $99.9\%$ of the energy on the OpInf subdomain and a nonlinear sFOM model (due to \eqref{eq:floating-condition}) with $n_F = 4,585$ DOFs by solving~\eqref{coupled_structure_linear}. In terms of subdomain selection, the aim is two-fold: the localized nonlinearity from the ice-ocean interaction in~\eqref{basal_melt} is active in the region close to the ocean, and, in parallel, the physical mechanisms in this region lead to transport-dominated dynamics for the ice thickness. The resulting sFOM and OpInf subdomains along with a narrow overlap region are shown in~\Cref{fig:ice1}. 

\begin{figure}[!ht]
    \centering
    \begin{minipage}{0.48\textwidth}
        \centering
        \includegraphics[width=0.99\textwidth]{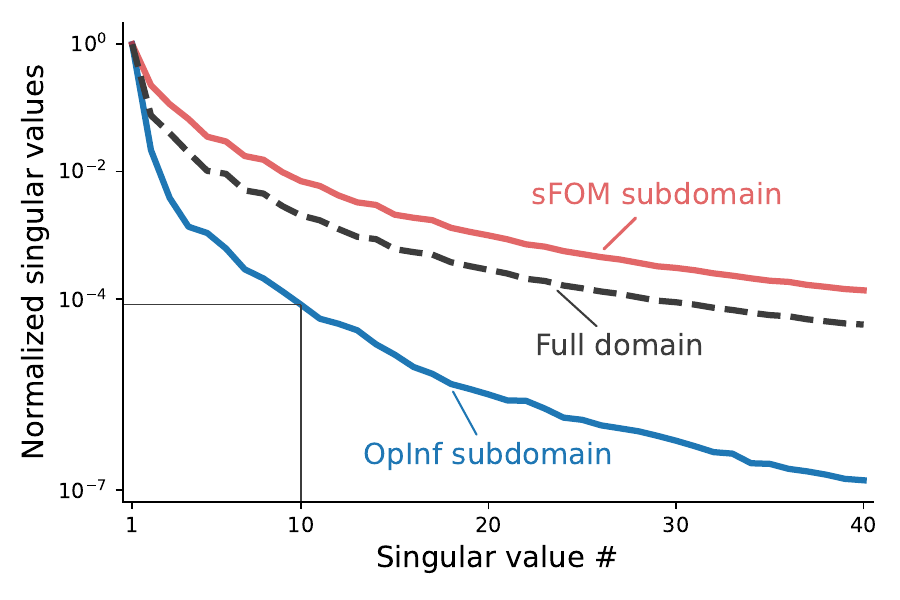}
\caption{Normalized singular value decay over training data. The data in the sFOM subdomain exhibit advective ice thickness dynamics and thus a slower singular value decay compared to the data in the OpInf subdomain.} 
\label{fig:svals_ice}
    \end{minipage}\hfill
    \begin{minipage}{0.48\textwidth}
        \centering
        \includegraphics[width=1.01\textwidth]{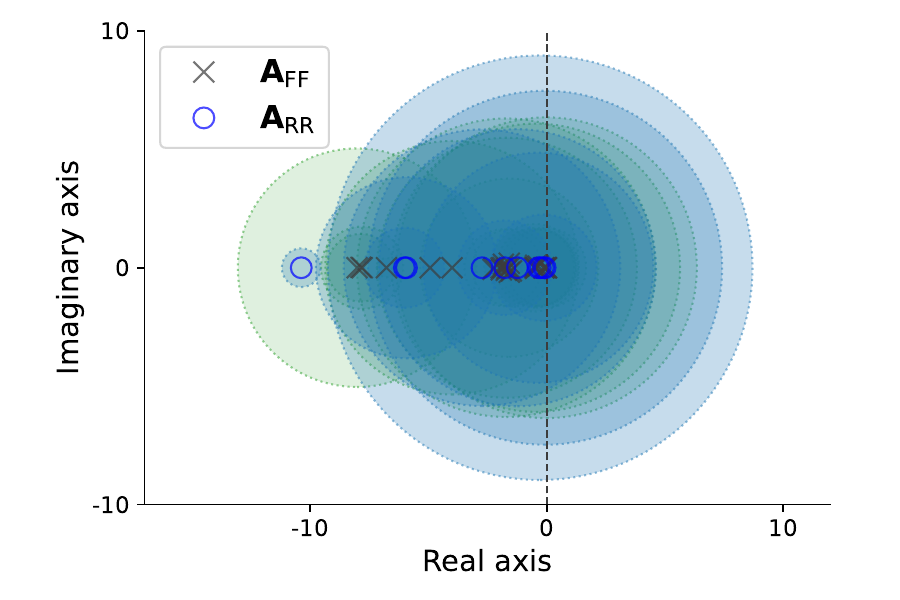}
\caption{Eigenvalues and Gershgorin disks for the inferred OpInf and sFOM linear operators $\mathbf{A}_{\text{FF}}$ and $\mathbf{A}_{\text{RR}}$ in \eqref{coupled_inference}, using the introduced regularization in \eqref{sfom_gersh} and the values in \Cref{tab:hyp_ice}. For $\mathbf{A}_{\text{FF}} \in \mathbb{R}^{4,458 \times 4,458}$ we plot 30, randomly selected eigenvalues. The eigenvalue of the inferred $\mathbf{A}_{\text{FF}}$ with the largest real part is $\lambda =0.051+0i$, while $\mathbf{A}_{\text{RR}}$ is stable.}
\label{fig:ice_eigs}
    \end{minipage}
\end{figure}

To train the parametric OpInf-sFOM, we use the first 15 years of simulation data for the extreme cases of $\alpha=\SI{0}{m/yr}$ and $\alpha=\SI{100}{m/yr}$ observed in \cite{rignot2013ice}. The sFOM subdomain is selected to include the area close to the ocean where the nonlinear input~\eqref{basal_melt} is active for the training data, as illustrated in \Cref{fig:ice1}. The corresponding singular value decay of both subdomains in \Cref{fig:svals_ice} supports this selection, since the singular value decay in the OpInf subdomain is faster than that in the sFOM subdomain. We choose the optimal regularization values for both models in~\eqref{sfom_gersh} by an L-curve criterion. For the inferred sFOM, we consider an L-curve optimization for both $\eta_1^{\text{sFOM}}$ and $\eta_2^{\text{sFOM}}$. This is done by constructing one L-curve per $\eta_2^{\text{sFOM}}$ value. The optimal pair$\left( \eta_1^{\text{sFOM}}, \eta_2^{\text{sFOM}} \right)$ then corresponds to the point (among all L-curves) that is closest to the origin of the normalized LS $L_2$ error against the normalized $L_2$ norm of the inferred $\bm{\beta}_i$ in \eqref{sfom_ls}. For OpInf, we optimize $\eta_1^{\text{OpInf}}$ and consider $\eta_2^{\text{OpInf}} = 2\times10^3\times\eta_1^{\text{OpInf}}$. The range of the regularization values for both OpInf and sFOM are summarized in \Cref{tab:hyp_ice}. \Cref{fig:ice_eigs} illustrates the effect of the Gershgorin regularization in~\eqref{sfom_gersh} to the OpInf and sFOM eigenvalues in~\eqref{coupled_inference} and showcases the robustness of the resulting parametric OpInf-sFOM.

\begin{table}[!ht]
  \caption{Coupled OpInf-sFOM inference hyperparameters for ice thickness predictions. The optimal regularization values in \eqref{sfom_gersh} are selected via an L-curve criterion \cite{Hansen2000} for the OpInf and the sFOM inference. For both OpInf and sFOM, $n_{\eta_1}$ logarithmically spaced values for $\eta_1$ were tested. For sFOM, $n_{\eta_2}$ logarithmically spaced values for $\eta_2$ were also tested.}
  \label{tab:hyp_ice}
  \centering
    \begin{tabular}{ccc|cc|cc}
    \hline\noalign{\smallskip}
      $r$ & $\eta_{1}^{\text{OpInf}}$ & $\eta_{2}^{\text{OpInf}}$ &$\eta_{1}^{\text{sFOM}}$ & $\eta_{2}^{\text{sFOM}}$ & $n_{\eta_1}$ & $n_{\eta_2}$\\
      \noalign{\smallskip}\hline\noalign{\medskip}
      10 & $\left[ 10^2, \; 10^7\right]$ & $2 \times 10^3 \times \eta_{1}^{\text{OpInf}}$ & $\left[ 10^2, \; 10^7\right]$ & $\left[ 10^{-2}, \; 10^{2}\right] \times \eta_{1}^{\text{sFOM}}$ & 20 & 8\\
      \noalign{\medskip}\hline
    \end{tabular}
  %\end{center}
\end{table}

We simulate the resulting coupled OpInf-sFOM of dimension $r+n_\text{F}=10+4,585=4,595$ using a 4th order Runge-Kutta integration scheme. The offline speedup compared to a global sFOM~\eqref{ratio_sfomg} is $2.27$, while a global OpInf would be faster by a factor of $\sim 13$~\eqref{ratio_opg}, irrespective of its predictive capabilities. We note that the coupled OpInf-sFOM efficiently handles the distributed nonlinearity due to~\eqref{eq:floating-condition} and the resulting transport-dominated dynamics in the sFOM subdomain, which would pose additional challenges to a global OpInf model. The theoretical online speedup against an explicitly available FOM code with the considerations made in~\eqref{online_cost_ratio} is $\sim 2.67$. In our case, simulating the coupled OpInf-sFOM on the same machine as the ISSM code for $21$ evenly distributed $\alpha$ values in $\left[0, \; 100 \right] \SI{}{m/yr}$, we observe an average online speedup factor of $7.89$. The discrepancy between the theoretically expected and practically observed speedup can be attributed to different software architectures, programming languages and  time stepping schemes used for the ISSM and the OpInf-sFOM, contrary to the considerations made for the theoretical analysis in~\eqref{online_cost_ratio}.

\begin{figure}[!ht]
\centering
\includegraphics[width=\textwidth]{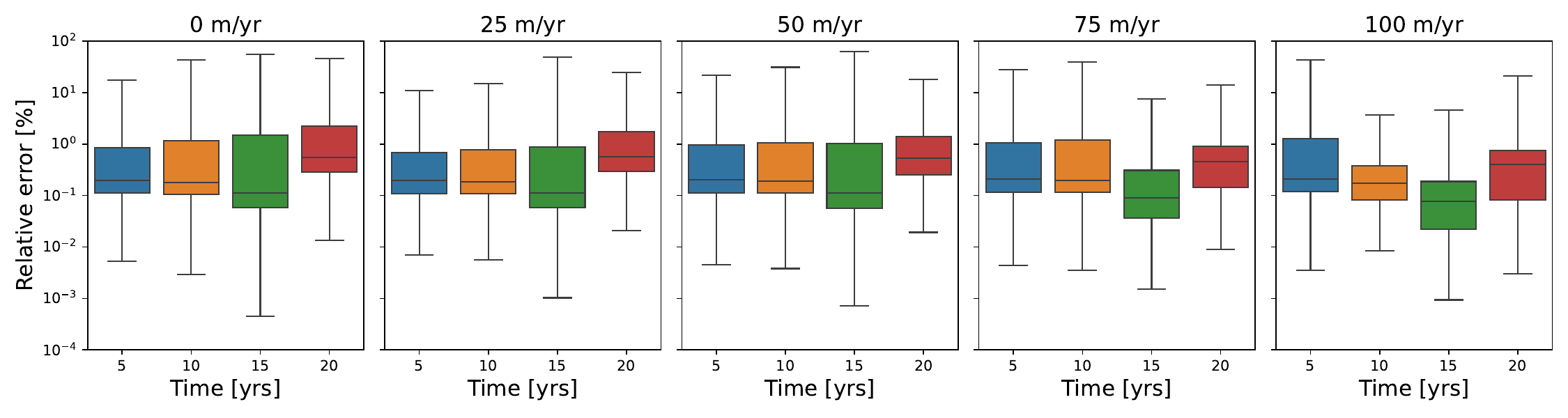}
\caption{State prediction error distribution over time, for selected $\alpha$ values. As time passes, the relative error median increases, with the boxplot distribution maximum reaching values of the order of $10 \%$. At each timestep, the error is normalized by the mean ice thickness across the spatial domain.}
\label{fig:ice_err_time}
\end{figure}

\begin{figure}[!t]
\centering
  \includegraphics[width=.7\textwidth]{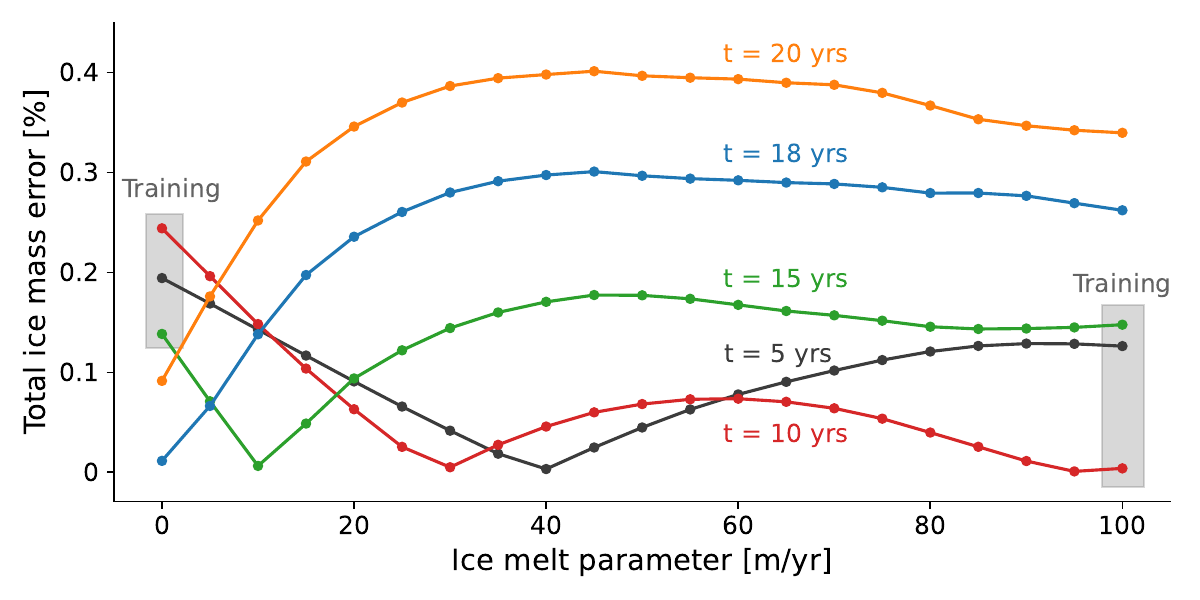}
\caption{Predicted total ice mass error. We observe low error values for all $\alpha$, with an increasing trend as time passes. The instances included in the training data are marked with a grey box.}
\label{fig:ice_mass}
\end{figure}

The distribution of the state prediction relative error for selected $\alpha$ values and over time $t$ is given in \Cref{fig:ice_err_time}. As a general trend, the median error increases over time for all parameters $\alpha$ up to values in the order of $1 \%$. By integrating the ice thickness over the complete domain, we also compute the predicted total ice mass for uniformly distributed $\alpha$ values. The corresponding error is plotted in \Cref{fig:ice_mass}. It exhibits an increasing trend over time for most of the parameter values, with an upper bound of $\sim 0.4 \%$. These results indicate that using the aforementioned training data only for the extreme values of $\alpha =$ \SI{0}{m/yr} and $\alpha =$ \SI{100}{m/yr} is sufficient for the parametric coupled OpInf-sFOM to successfully predict the system dynamics for any intermediate $\alpha$ value. 

\begin{figure}[!t]
\centering
  \includegraphics[width=.7\textwidth]{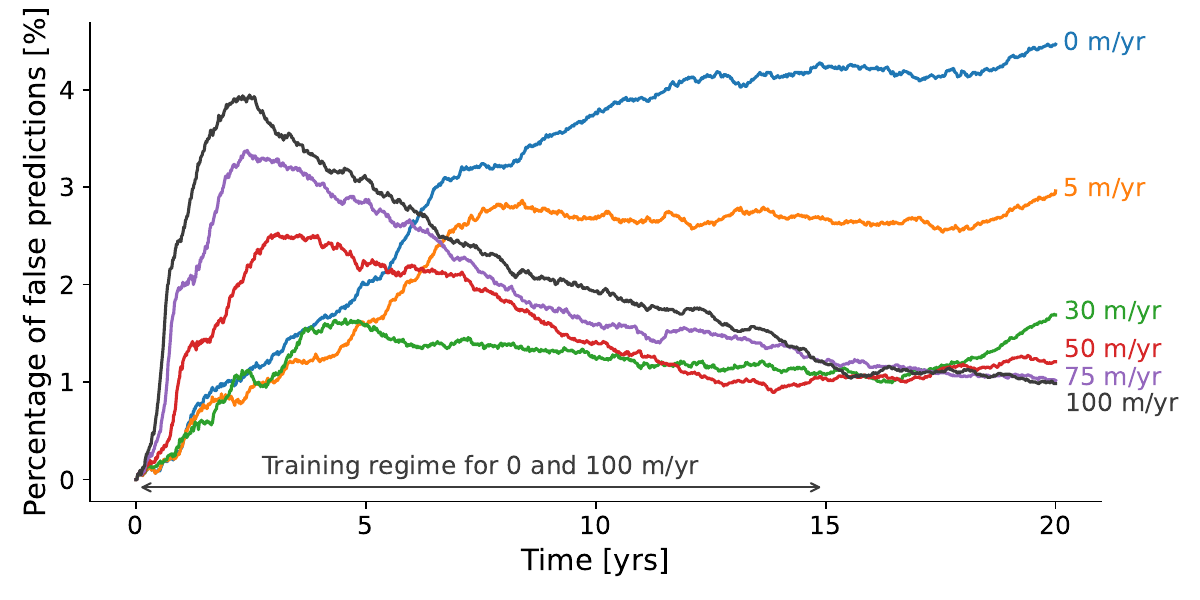}
\caption{Percentage of total DOFs falsely predicted to be floating or grounded. The largest error $\left( \sim 4.5 \% \right)$ is observed for $\alpha=$ \SI{0}{m/yr}, with no evident increasing trend beyond the $15$ years of training data. The error decreases for higher $\alpha$ values.}
\label{fig:false_floating}
\end{figure}

\begin{figure}[!ht]
\centering
\includegraphics[trim=130 100 95 80, clip, width=.95\textwidth]{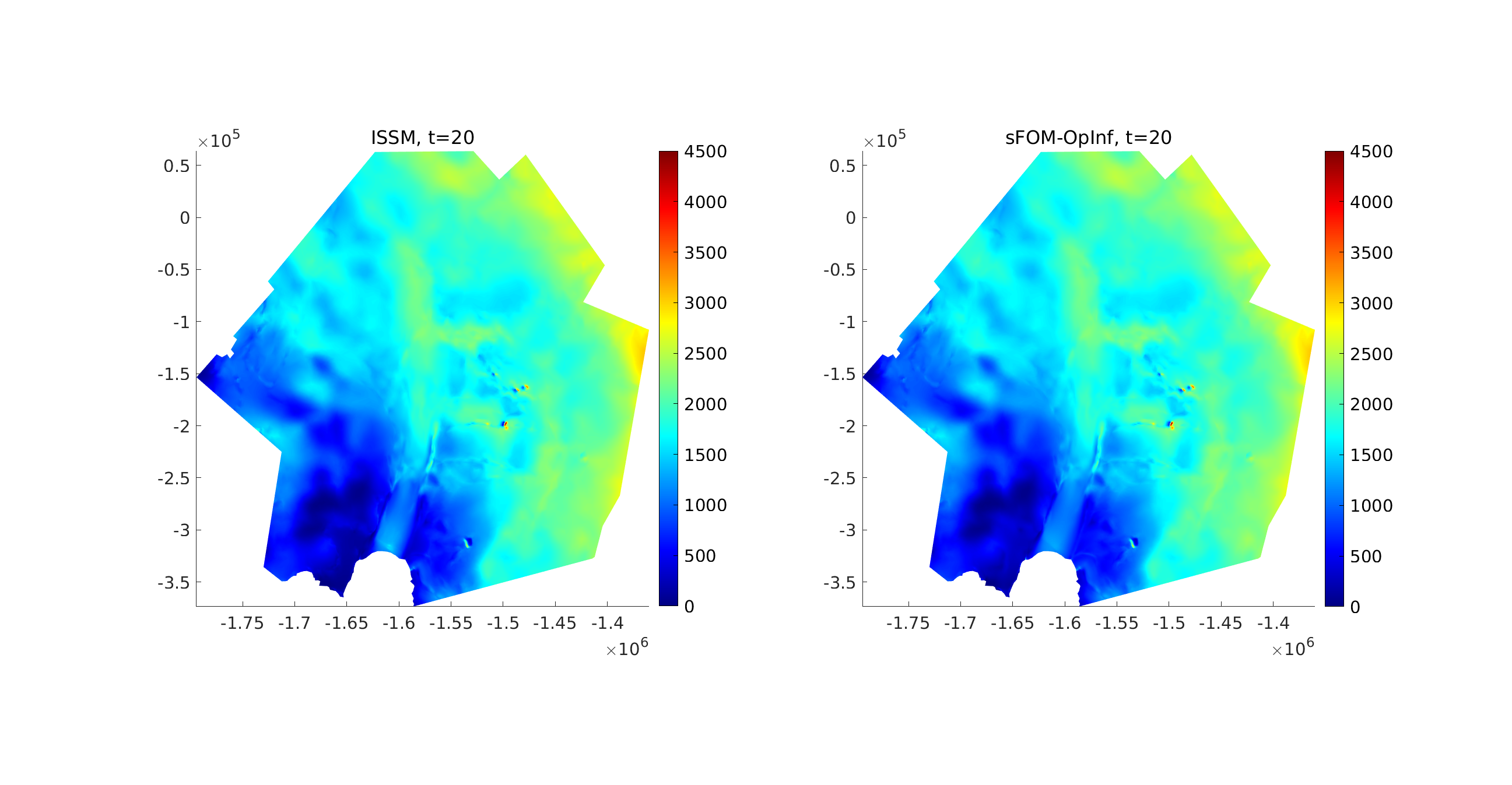}
\vfill   \includegraphics[trim=130 80 95 50, clip, width=.95\textwidth]{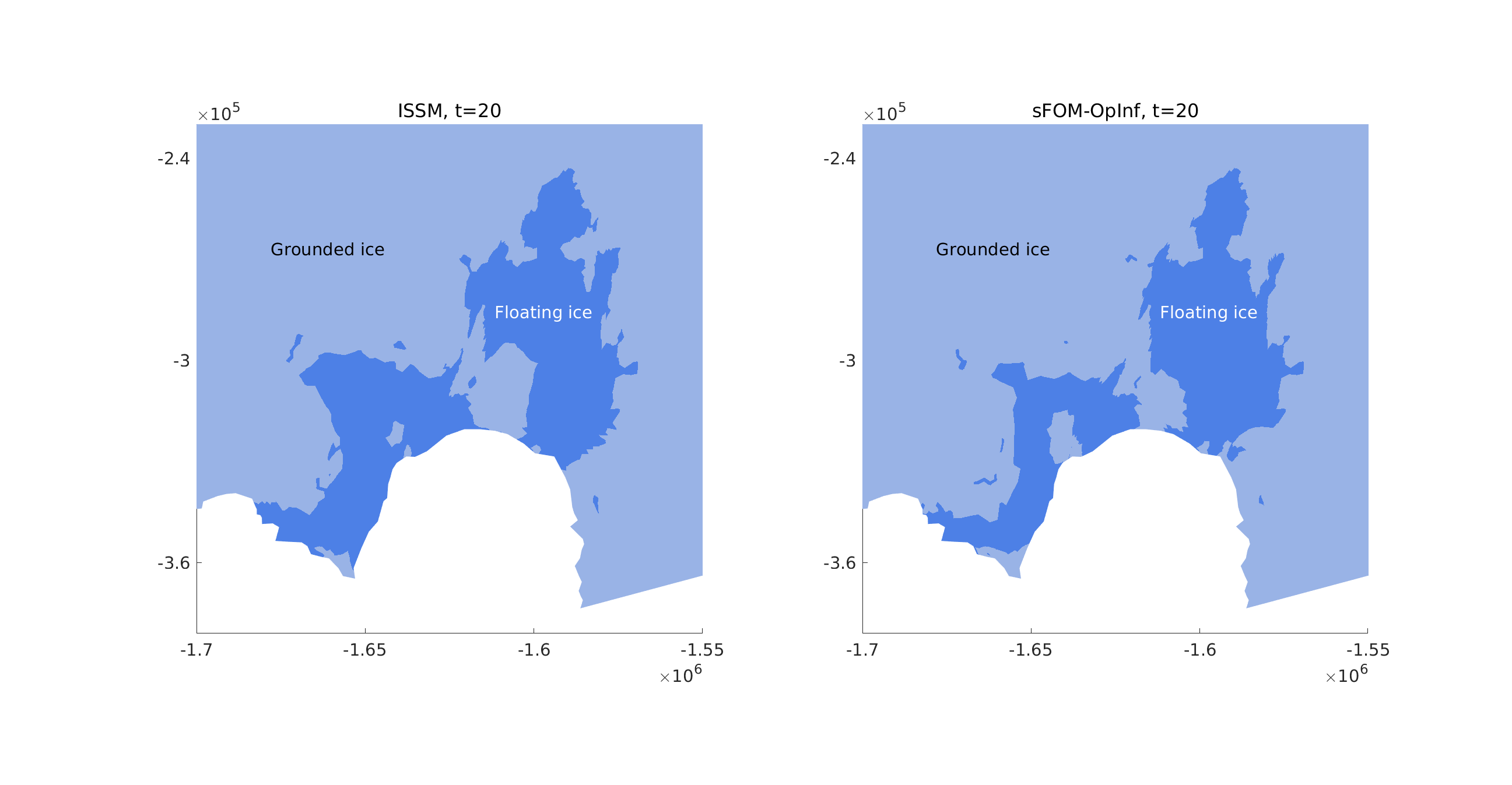}
\caption{OpInf-sFOM predictions for $\alpha=$ \SI{0}{m/yr}. Upper row: state predictions for the coupled OpInf-sFOM and comparison to the ISSM code data. Lower row: predicted floating ice region, with a close-up towards the sFOM subdomain. The predictions are quantitatively matching the data, with $4.5\%$ of the total DOFs being falsely classified as floating/grounded.}
\label{fig:0_states}
\end{figure}

\begin{figure}[!htb]
\centering
  \includegraphics[trim=130 100 95 80, clip, width=.95\textwidth]{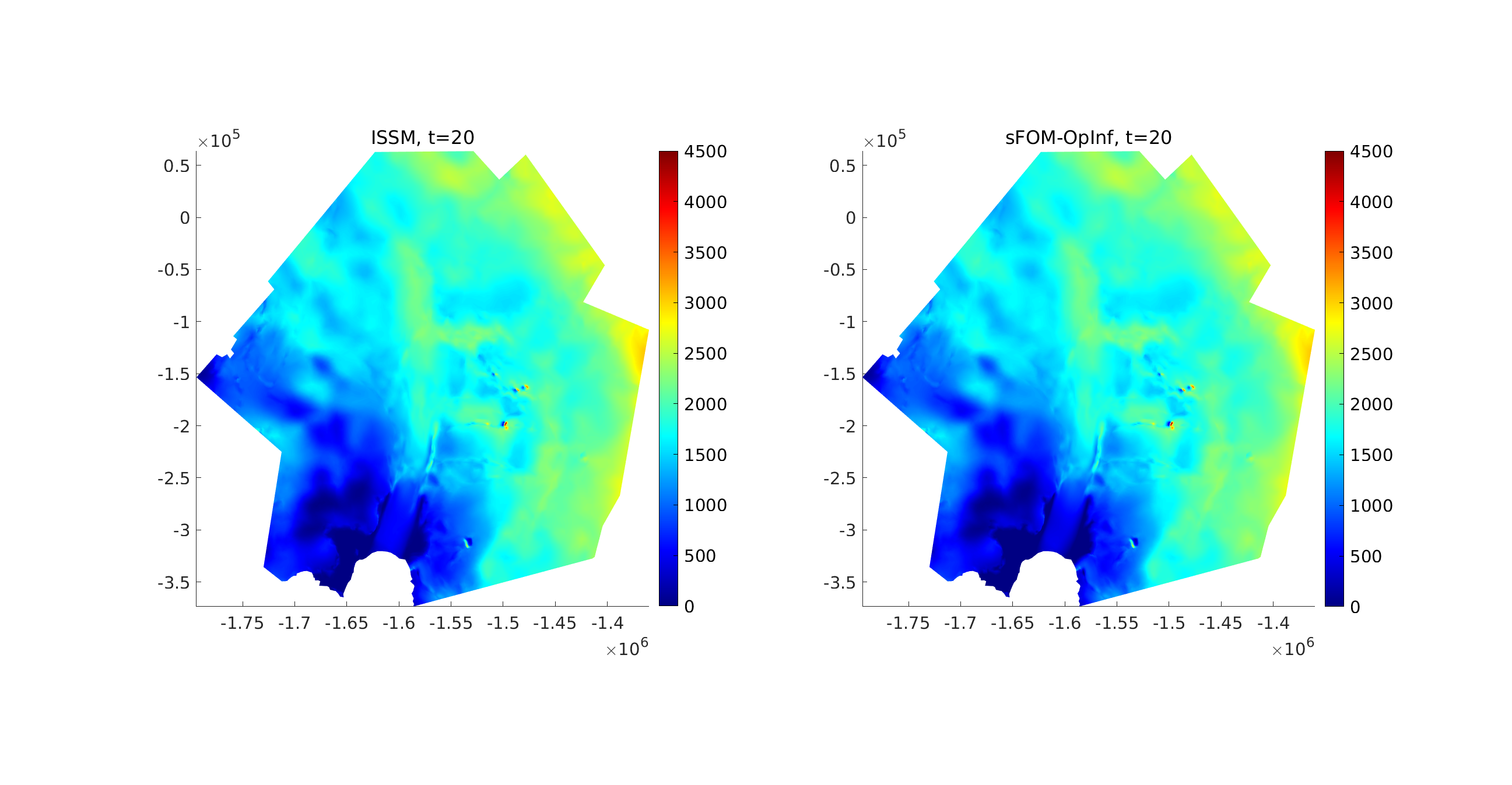}
  \includegraphics[trim=130 80 95 50, clip, width=.95\textwidth]{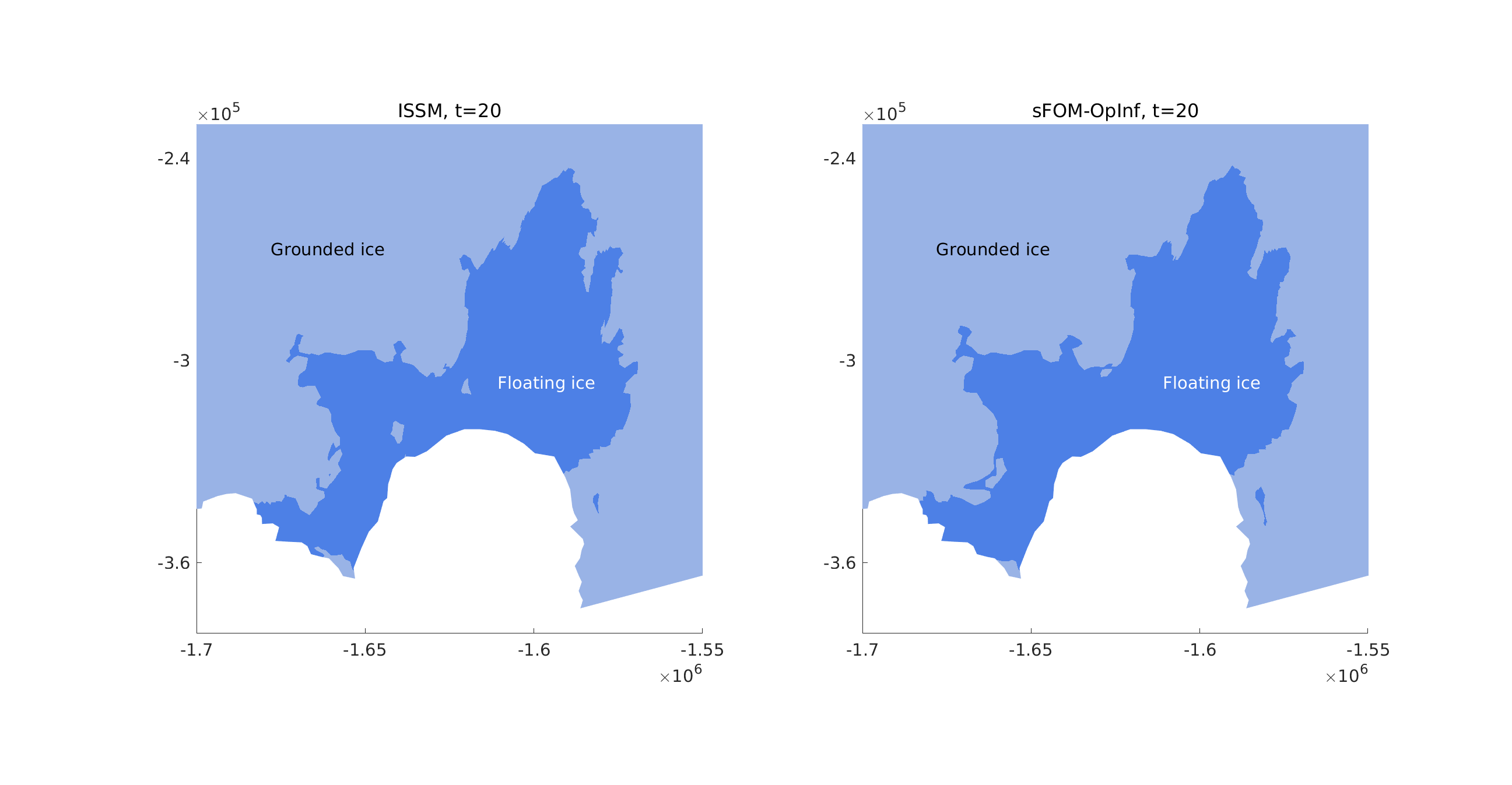}
\caption{OpInf-sFOM predictions for $\alpha=$ \SI{50}{m/yr}. Upper row: state predictions for the coupled OpInf-sFOM and comparison to the ISSM code data. Since the largest part of the DOFs where input \eqref{basal_melt} is applied are in the sFOM subdomain, the OpInf has no parametric dependence. Lower row: predicted floating ice region, with a close-up towards the sFOM subdomain. The predictions are in good accordance with the data, with $1\%$ of the total DOFs being falsely classified as floating/grounded.}
\label{fig:50_states}
\end{figure}

An important metric in ice dynamics is the prediction of the floating ice regions, following the criterion in \eqref{eq:floating-condition}. For this reason, we compute the percentage of DOFs that the parametric OpInf-sFOM classifies falsely to be either floating or grounded ice. \Cref{fig:false_floating} presents the evolution of this false prediction in time. We observe that for $\alpha=$ \SI{0}{m/yr} we have the highest error of approximately $4.5 \%$. To make these values tangible, we present the comparison of the state predictions to the CFD data at $t=20$ years and the corresponding comparison for floating ice. \Cref{fig:0_states} corresponds to $\alpha=$ \SI{0}{m/yr}, for which the floating/grounded DOFs error is the highest in \Cref{fig:false_floating}, while \Cref{fig:50_states} corresponds to $\alpha=$ \SI{50}{m/yr}, for which the floating/grounded DOFs error is below $1.5 \%$ at $t=20$. We conclude that the coupled OpInf-sFOM, trained only with the first 15 years of simulation data for $\alpha=\SI{0}{m/yr}$ and $\alpha=\SI{100}{m/yr}$, serves as an accurate surrogate model for ice thickness dynamics predictions in the parameter range of $\alpha\in\left[0,\; 100\right]\SI{}{m/yr}$.

\section{Conclusions and future perspectives}
\label{sec:conclusion}

In this paper we presented a non-intrusive coupled ROM/FOM approach designed for the inference of dynamical systems with spatially localized slow singular value decay, building upon the methods of Operator Inference~\cite{peherstorfer2016data} and sparse FOM inference~\cite{SCHUMANN2023}. We presented the two methods, indicated their advantages and limitations and formulated a coupled OpInf-sFOM with the use of domain decomposition. We also discussed the usage of the spectral gap indicator for validating the OpInf and sFOM subdomains selection and examined the issue of solution smoothness over the OpInf-sFOM interface. The subsequent parametric analysis for the online and offline speedup by the coupled OpInf-sFOM quantified the impact of the reducibility in the OpInf subdomain and the sFOM subdomain size to the computational efficiency of the method.
Moreover, we introduced the stability-promoting Gershgorin regularization and provided a closed-form solution for the inference LS problem in both OpInf and sFOM. This advancement was crucial for inferring robust non-intrusive dynamical systems. The potential of the coupled OpInf-sFOM formulation with Gershgorin regularization was exemplified for two numerical test cases. The 1D Burgers' example illustrated the potential for predictions beyond the training regime for transport-dominated dynamics. The 2D ice thickness dynamics model for the Pine Island Glacier further showcased the accuracy and speedup for parametric predictions in a realistic application with a localized nonlinearity. The inferred OpInf-sFOM was found to be $7.89$ times faster than the FOM code, while achieving an average error of the order of $1 \%$ for the range of examined parameters.
 
The current formulation can be extended in several directions, towards a comprehensive framework encompassing physics-informed non-intrusive modeling and data assimilation. First, we consider an extension of the presented coupled approach towards an online, adaptive decomposition of the domain to ROM and FOM subdomains. This approach would enable predictions for systems with dynamically evolving regions of slow singular value decay, while requiring less prior knowledge about the localization of such features. Moreover, an online adaptive ROM-FOM decomposition would allow for bigger speedups in systems with fixed regions exhibiting slow singular value decay, while maintaining the flexibility of localized inference at the full-order level, where and when necessary. Inference at the full-order level could also be leveraged for the incorporation of measurement data from sensors, in the direction of data assimilation. Both aforementioned directions require a further investigation of data-driven indicators and corresponding algorithms for domain decomposition. Finally, testing the OpInf-sFOM approach for the inference of systems undergoing bifurcations is an interesting direction for future research. For such systems, the full-order level inference on the sFOM subdomain could allow for an accurate prediction of local, bifurcating dynamics, while retaining a considerable speedup via the inferred OpInf model on the corresponding subdomain. 

\section*{Data $\&$ Code Availability}

The simulation data and Python code for the application of the coupled OpInf-sFOM methodology to the 1D Burgers' equation test case in \Cref{subs:1d_burg} are available in \url{https://github.com/lgkimisis/opinf_sfom}.

\section*{Acknowledgements}

Leonidas Gkimisis, Thomas Richter and Peter Benner acknowledge support from the Deutsche Forschungsgemeinschaft (DFG, German Research Foundation) - 314838170, RTG 2297 MathCoRe. 
The authors from the University of Texas at Austin acknowledge support from the Department of Energy grants DE-SC0021239 and DE-SC002317, the Air Force Office of Scientific Research grant FA9550-21-1-0084, and the NASA University Leadership Initiative under Cooperative Agreement 80NSSC21M0071.

\begin{appendix}
\section{OpInf-sFOM interface position}
\label{appendix:app}

\begin{figure}[!htb]
\centering
  \includegraphics[width=.75\textwidth]{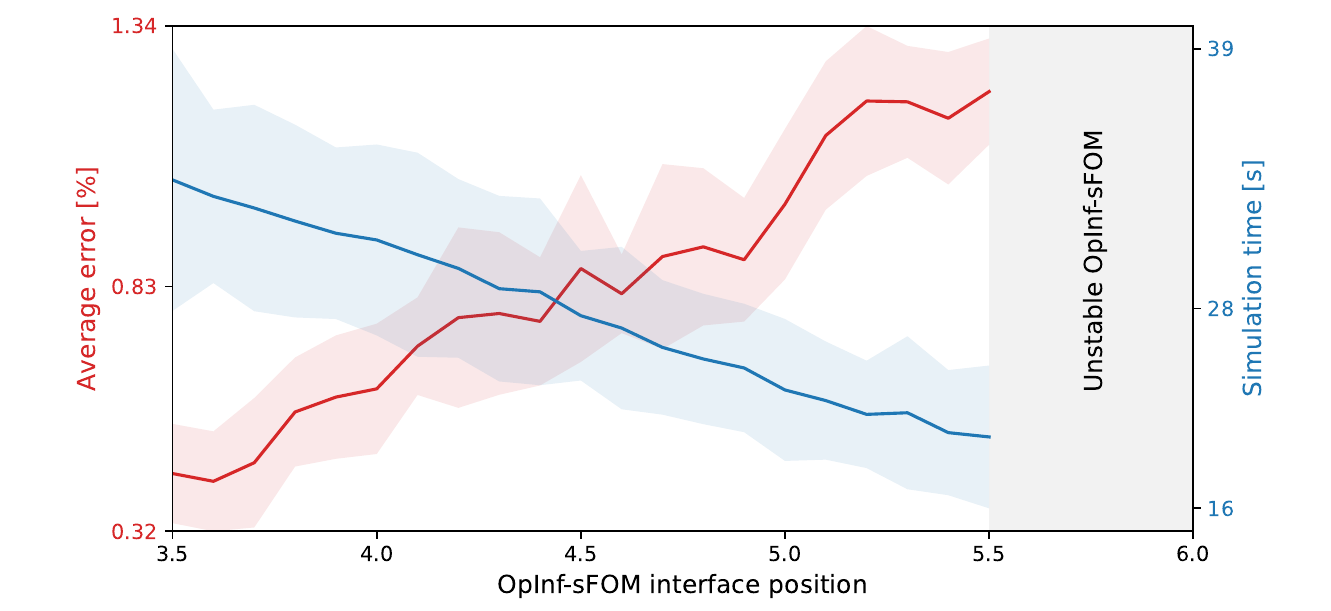}
\caption{Effect of OpInf-sFOM interface position on coupled, non-intrusive predictions for 1D Burgers'. As the sFOM region shrinks, the OpInf-sFOM simulation time decreases and the average error over testing time and space increases. However, unstable OpInf-sFOM models can arise when a significant portion of the transport-dominated region in \Cref{fig:1d_burg} is attributed to OpInf. The mean and $90\%$ confidence interval of the results over $10$ simulations for each value of $a$ are given.}
\label{fig:burg_param}
\end{figure}

We perform a parametric study with respect to the position of the OpInf-sFOM interface for the 1D Burgers' test case in \Cref{fig:1d_burg}. In particular, we consider an OpInf model for $x \in [0, \;a ]$ and an sFOM model for $x \in [a, \; 10]$. We consider different values of $a$ in $a \in [3.5, 5.5]$ with a step of $0.01$ and learn coupled OpInf-sFOMs, while keeping all other training and inference parameters in \Cref{subs:1d_burg} constant. We thus gain insight on the effect of the domain partitioning on the properties of the coupled, inferred model. The parametric results on the OpInf-sFOM simulation time and the average error in space over the testing time are given in \Cref{fig:burg_param}. All models were trained and run $10$ times to account for the variability in simulation time and error (due to the augmentation of \eqref{sfom_ls} with data from 5 random spatial points). As the position of the interface $a$ increases, the average error of the OpInf-sFOM model increases, since a larger subdomain is attributed to OpInf. In parallel, this leads to a decrease in the required computational time for the OpInf-sFOM simulation. However, increasing the interface position beyond $a=5.5$ leads to unstable models (for the tested regularization values). Since the OpInf subdomain includes a region exhibiting transport-dominated dynamics (see \Cref{fig:1d_burg}), the $r=10$-dimensional linear basis leads to high projection errors and thus unsuccessful inference of the coupled dynamics.

\end{appendix}

\bibliographystyle{abbrvurl}
\bibliography{sfom_opinf.bib}

\end{document}